\documentclass[12pt]{amsart}

\usepackage{tgpagella}
\usepackage{euler}
\usepackage[T1]{fontenc}
\usepackage{amsmath, amssymb, mathabx}
\usepackage[hidelinks]{hyperref}
\usepackage[english]{babel}
\usepackage{mathrsfs}
\usepackage{eucal}
\usepackage[all]{xy}
\usepackage{tikz}

\newtheorem{thm}{Theorem}[section]
\newtheorem*{thm*}{Theorem}
\newtheorem{lem}[thm]{Lemma}

\newtheorem{prop}[thm]{Proposition}
\newtheorem*{prop*}{Proposition}
\newtheorem{conj}[thm]{Conjecture}
\newtheorem{cor}[thm]{Corollary}
\newtheorem*{cor*}{Corollary}

\theoremstyle{definition}
\newtheorem{defn}[thm]{Definition}
\newtheorem*{defn*}{Definition}

\newtheorem{remark}[thm]{Remark}

\newtheorem{example}[thm]{Example}

\newtheorem*{question*}{Question}
\newtheorem*{Pquestion*}{Popa's question}

\newtheorem*{conv*}{Convention}

\def\bb{\mathbb}

\def\bb{\mathbb}

\def\cal{\mathcal}
\newcommand{\cA}{\mathcal{A}}
\newcommand{\cB}{\mathcal{B}}

\makeatletter

\def\dotminussym#1#2{%
  \setbox0=\hbox{$\m@th#1-$}%
  \kern.5\wd0%
  \hbox to 0pt{\hss\hbox{$\m@th#1-$}\hss}%
  \raise.6\ht0\hbox to 0pt{\hss$\m@th#1.$\hss}%
  \kern.5\wd0}

\newcommand\cL{{\cal L}}
\newcommand\cC{{\cal C}}
\newcommand\bN{{\mathbb N}}
\newcommand\bR{{\mathbb R}}
\newcommand\bZ{{\mathbb Z}}

\allowdisplaybreaks
\makeatletter
\def\l@subsection{\@tocline{2}{0pt}{2.5pc}{5pc}{}}
\def\l@subsubsection{\@tocline{2}{0pt}{5pc}{7.5pc}{}}
\makeatother

\begin{document}


\title{Subsets of virtually nilpotent groups with the SBM property}

\author{Ryan Burkhart and Isaac Goldbring}
\address{Department of Mathematics\\University of California, Irvine, 340 Rowland Hall (Bldg.\# 400),
Irvine, CA 92697-3875}
\email{isaac@math.uci.edu}
\urladdr{http://www.math.uci.edu/~isaac}
\thanks{The results of this paper are taken from Burkhart's doctoral dissertation, written under Goldbring's supervision.}
\thanks{Goldbring was partially supported by NSF grant DMS-2054477.}

\maketitle

\begin{abstract}
We extend Leth's notion of subsets of the integers satisfying the Standard interval measure (SIM) property to the class of virtually nilpotent groups and name the corresponding property the Standard ball measure (SBM) property.  In order to do this, we define a natural measure on closed balls in asymptotic cones associated to such groups and show that this measure satisfies the Lebesgue density theorem.  We then prove analogs of various properties known to hold for SIM sets in this broader context, occasionally assuming extra properties of the group, such as the small spheres property and the small gaps property.
\end{abstract}

\tableofcontents

\section{Introduction}

In the 1980s, Steven Leth introduced a property of subsets of the natural numbers that he called the \textbf{SIM} (short for \textbf{standard interval-measure}) property \cite{LethSIM}.  The property is motivated by ideas from nonstandard analysis and, roughly speaking, demands that there be a connection between the \emph{internal} notion of small gap sizes on hyperfinite intervals and the \emph{external} notion of the (normalized) standard part of the set having large Lebesgue measure.  More precisely, given an infinite, hyperfinite interval $I=[y,z]\subseteq {}^*\bb N$, one can consider the map $\operatorname{st}_I:I\to [0,1]$ given by $\operatorname{st}_I(x):=\operatorname{st}\left(\frac{x-y}{z-y}\right)$.  It can be shown that if an internal subset $A$ of $I$ is such that $\operatorname{st}_I(A)$ has large Lebesgue measure, then all gaps of $A$ on $I$ have small size (that is, small internal cardinality) compared to the size of $I$.  The SIM property is motivated by asking that the converse relation hold, uniformly over all infinite, hyperfinite subintervals of $I$.

In the same paper, Leth showed that subsets of the natural numbers with the SIM property satisfy the conclusion of a theorem of Steward and Tijdeman \cite{StewartTijdeman} known to hold for sets with positive Banach density:  if $A_1,\ldots,A_n\subseteq \mathbb{N}$ have the SIM property, then $D(A_1)\cap \cdots\cap D(A_n)$ is syndetic, where, for any subset $A\subseteq \bb N$, $D(A)$ is the set of those $n\in \bb N$ which may be written as $a_i-a_j$ for infinitely many pairs $(a_i,a_j)\in A^2$.

In joint work of Leth and the second author, further properties of SIM sets were established.  For example, Leth had shown that being a SIM set is not simply a notion of size and it became apparent that it would be sensible to consider the class of \textbf{supra-SIM} sets, that is, the class of sets containing SIM sets.  In \cite{LGsupraSIM}, it was shown that the class of supra-SIM sets is a partition regular class, whence is combinatorially natural.  Moreover, they were able to prove a version of Jin's theorem for sumsets of SIM sets, namely the sum of any two SIM sets is piecewise syndetic \cite{Jin} (Jin's theorem had the same conclusion but with $A$ and $B$ assumed to have positive Banach density), as well as a version of Nathanson's theorem from \cite{Nathanson}, namely if $A$ is SIM, then for any $n\in \bb N$, there are $B,C\subseteq \bb N$ with $B$ infinite and $|C|=n$ for which $B+C\subseteq A$.  (Nathanson's result had the same conclusion under the assumption that $A$ had positive Banach density.)  Nathanson's theorem was partial progress on a sumset conjecture of Erdos, which asked if $A$ has positive lower density, must $A$ contain the sum of two infinite sets?  This conjecture was recently resolved in the affirmative, under the weaker assumption of positive Banach density, by Donaldson, Moreira, and Richter \cite{SumsetConj}.  Whether or not the theorem holds under the assumption that $A$ is simply a SIM set is an interesting open question.

It is relatively routine to generalize the definition of SIM set from subsets of $\bb N$ to subsets of $\bb Z$.  By replacing intervals by cubes, one can easily generalize the definition further to subsets of $\bb Z^n$ and thus to subsets of any finitely generated abelian group.  However, generalizing the definition beyond the finitely generated abelian case is not entirely clear as many of the arguments rely heavily on the ``geometry'' of intervals and cubes.

One idea of how to expand the class of groups to which this concept is meaningful is to replace cubes by balls in the Cayley graph (with respect to some given finite set of generators).  The natural quotient map from an infinite, hyperfinite ball to the \textbf{asymptotic cone} of the group with respect to the given set of generators and some infinite hypernatural number is then a suitable replacement for the map $\operatorname{st}_I$ considered above.  However, many of the results mentioned above rely heavily on the Lebesgue density theorem for Lebesgue measure on $[0,1]$.  Since the Lebesgue measure is the pushforward of the Loeb measure associated to the hyperfinite counting measure on $I$ via $\operatorname{st}_I$, it is thus natural in our extended setting to ask that the Lebesgue density theorem hold for the pushforward of the Loeb measure associated to the hyperfinite counting measure on infinite, hyperfinite balls in the Cayley graph.  We show that, assuming that the group has \textbf{polynomial growth}, this pushforward measure does indeed satisfy the Lebesgue density theorem.  By Gromov's celebrated theorem \cite{GromovThm}, the assumption of polynomial growth is equivalent to the assumption that the group is virtually nilpotent.

In this paper, we show how to extend the definition of SIM sets to the class of finitely generated virtually nilpotent groups, obtaining a class of sets that we call \textbf{SBM sets}, where SBM is short for the \textbf{standard ball measure} property.  However, to obtain results that parallel the results discussed above, we sometimes need to assume that our groups satisfy a property we call the \textbf{small sphere} property, a property known to hold in all nilpotent groups of class 2 and conjectured to hold for all nilpotent groups.  For some of our results, a third property, which we call the \textbf{small gap} property, is also assumed.  With these assumptions in hand, we prove direct analogs of the aformentioned results known to hold for SIM sets (with the SBM analog of Jin's theorem being the only exception).

While it is satisfying to extend the class of SIM sets to a fairly large class of nonabelian groups, it would be desirable to extend the notion even further.  In particular, it would be interesting to see if this notion could be suitably extended to the class of amenable groups, perhaps using well-behaved Folner sequences as replacements for balls in the Cayley graph.  Ideally one could even extend the notion past amenable groups.  For example, how to extend the notion to finitely generated free groups is also not clear to us at present.

We assume the reader is familiar with the basics of nonstandard analysis.  The monograph \cite{GoldbringBook} contains a complete introduction and is aimed at an audience interested in combinatorial applications.  Many of the proofs appearing later in this paper are based on the arguments given in Chapter 16 of that monograph.  Nevertheless, some basic notation and information on certain nonstandard constructions, such as hyperfinite products in groups and the Loeb measure construction, are given in the next section.

\section{Preliminaries}

\subsection{Some nonstandard analysis}

In this subsection, we only recall a few facts from nonstandard analysis that we use throughout the paper.  We first remind the reader that, for any set $X$, letting $\mathcal{P}(X)$ denote the powerset of $X$, we can naturally identify the nonstandard extension $^*\cal P(X)$ of $\cal P(X)$ with a subset of $\cal P(X^*)$; in this way, elements of $^*\cal P(X)$ will be called \textbf{internal subsets of $^*X$}; subsets of $^*X$ that are not internal are called external.  In a similar way, if $\cal P_f(X)$ denotes the finite subsets of $X$, then elements of $^*\cal P_f(X)$ are special internal subsets of $^*X$ called \textbf{hyperfinite} subsets of $^*X$.  By transfer, for each hyperfinite subset $A$ of $^*X$, there is a unique element $N$ of $^*\bb N$ for which there is an internal bijection $f:A\to [1,N]:=\{M\in {}^*\bb N \ : \ 1\leq M\leq N\}$; we refer to this $N$ as the \textbf{internal cardinality} of $A$, denoted by $|A|$.  Note also that, by transfer, an internal subset of a hyperfinite set is once again hyperfinite.

Given any hyperfinite set $A$, there is a natural finitely-additive probability measure $\mu$ on the algebra (but not $\sigma$-algebra!) of internal subsets of $A$ given by $\mu(B):=\operatorname{st}\left(\frac{|B|}{|A|}\right)$.  It is a standard fact that the hypotheses of the Carath\'{e}odory extension theorem hold in this context, whence we get an extension of $\mu$, denoted $\mu_L$ and called the \textbf{Loeb measure}, defined on a canonical $\sigma$-algebra of subsets of $A$ containing the internal subsets, called the $\sigma$-algebra of \textbf{Loeb measurable sets}.  An example to keep in mind is the following:  if $N\in {}^*\bb N\setminus \bb N$ and $A=\{0,\frac{1}{N},\ldots,\frac{N-1}{N},1\}$ and $\operatorname{st}:A\to [0,1]$ is the usual standard part map, then the pushforward of the Loeb measure space on $A$ via $\operatorname{st}$ is the usual Lebesgue measure on $[0,1]$. 

We also remind the reader of the notion of hyperfinite products in groups.  Suppose that $G$ is an arbitrary group.  Let $G_{fs}$ denote the set of finite sequences from $G$ and let $p:G_{fs}\to G$ denote the product map $p((g_1,\ldots,g_n)):=g_1\cdots g_n$.  We can thus consider the nonstandard extension of $^*G_{fs}$, which is the set of hyperfinite sequences from $^*G$; we suggestively write elements of $^*G_{fs}$ as $(g_1,\ldots,g_N)$ for some element $N\in {}^*\bb N$.  Similarly, we consider the nonstandard extension $p:{}^*G_{fs}\to {}^*G$ of $p$ and suggestively write $g_1\cdots g_N$ instead of $p((g_1,\ldots,g_N))$.

\subsection{Growth rates for finitely generated groups}
\

From now on, throughout this paper, $G$ will always denote a finitely generated group with finite generating set $S$.  For simplicity, we will always assume that our generating set $S$ is symmetric, that is, contains the identity and is closed under inverses.

The \textbf{word metric} on $G$ with respect to $S$ is the metric $$d_S(x,y) = \min\{n \in \bN : x^{-1}y = s_1s_2\cdot \cdot \cdot s_n \text{ for some } s_1,...,s_n \in S\}.$$

Note that $d_S$ is left-invariant, that is, for all $x,y,g \in G$, $d_S(x,y) = d_S(gx,gy)$. Also notice that, for $\bZ$ with the standard symmetric generating set $S = \{1, -1\},$ we have that $d_S$ is simply the usual notion of distance in $\bZ$.  

For any $g \in G$ and any $n \in \bN$, define the \textbf{closed ball of radius $n$ and center $g$} to be $B(g,n) = \{h \in G : |g^{-1}h| \leq n\}$, where we also allow $B(g,0) = \{g\}$. In what follows, we will often omit the word ``closed'' and simply refer to this as a ``ball'. The \textbf{open ball of radius $n$ with center $g$} is then defined as $U(g,n) = \{h \in G : |g^{-1}h| < n\} = B(g,n-1)$, allowing us to define the \textbf{sphere of radius $n$ with center $g$} by $S(g,n)=B(g,n) \backslash B(g,n-1)$.  Note that we can (and will) also allow real (instead of natural) number radii, noting that spheres with non-integer radii are empty.

For $g \in G$, set the \textbf{word length} of $g$ with respect to $S$ to be $|g|_S = d_S(g,e)$, where $e$ denotes the identity of the group. For $n \in \bN$, set $$\alpha_n =\alpha_{n,S}= |\{g \in G : |g|_S \leq n\}|=|B(g,e)|.$$  Note that we have the trivial upper bound $\alpha_n \leq |S|^n$. However, in many groups, $\alpha_n$ will grow much more slowly than this exponential upper bound. For example, if $G$ is abelian, then it is readily verified that $\alpha_n < (n+1)^{|S|}$ for large enough $n$.  This is an example of the following phenomenon:

\begin{defn}
We say that $G$ has \textbf{polynomial growth} if there is $d\in \bb N$ and $C>0$ such that $\alpha_n\leq Cn^d$ for all $n\in \bb N$. The \textbf{degree} of polynomial growth of a group is defined to be the minimum degree $d$ for which the above inequality holds for some positive constant $C$.
\end{defn}

Since any two word metrics on $G$ are bi-Lipschitz equivalent, the notion of polynomial growth is indeed independent of the choice of generating set.

In many cases, it will be useful to have both upper and lower bounds for $\alpha_n$. 

\begin{thm} \cite{PolyDegree}\label{witness}
If $G$ has polynomial growth with degree $d$, then there exists a constant $c \in \bN$ such that $\frac{1}{c}n^d \leq \alpha_n \leq cn^d$ for all $n \in \bN$.
\end{thm}







It was shown by Wolf in \cite{WolfPoly} that all virtually nilpotent groups are of polynomial growth. Thirteen years later, Gromov proved the converse in \cite{GromovThm}, thereby completely characterizing groups with polynomial growth. 

\begin{thm}[\textbf{Gromov \cite{GromovThm}}]
The finitely generated groups with polynomial growth are precisely the finitely generated groups which are virtually nilpotent.
\end{thm}

By Gromov's theorem, in the sequel, we use the terms polynomial growth and virtually nilpotent interchangeably. 

We extend the above notions to the nonstandard setting as follows.  The metric $d_S$ on $G$ extends to an internal metric, also denoted $d_S$, on $^*G$, which, by transfer, satisfies, for all $x,y\in ^*G$: 
$$d_S(x,y)=\min\{N \in \text{ } ^*\bN : x^{-1}y = s_1s_2\cdots s_N \text{ for some } s_1,...,s_N \in S\}.$$
Here, $(s_1,s_2,\ldots, s_N)$ ranges over hyperfinite sequences from $S$ and $s_1s_2\cdots s_N$ is the corresponding hyperfinite product, as introduced in the previous subsection.  Similarly, for any $g\in {}^*G$ and any $R\in {}^*\bb R^{>0}$, we have notions of internally closed balls $B(g,R)$, internally open balls $U(g,R)$, and internal spheres $S(g,R)$.

We note for future reference one nonstandard consequence of the polynomial growth assumption.

\begin{lem}
Suppose that $G$ has polynomial growth.  Let $M \leq N$ be infinite elements of $^*\bN$. Then $\frac{\alpha_M}{\alpha_N} \approx 0 \iff \frac{M}{N} \approx 0$. 

\end{lem}

\begin{proof}
Let $d$ be the degree of polynomial growth of $G$, as witnessed by the constant $c$. First assume that $\frac{M}{N} \approx 0$. Then since $\alpha_M \leq cM^d$ and $\alpha_N \geq \frac{1}{c}N^d$, we get that $\frac{\alpha_M}{\alpha_N} \leq \frac{c^2M^d}{N^d} = c^2(\frac{M}{N})^d$, which is infinitesimal. 

For the other direction, assume $\frac{\alpha_M}{\alpha_N} \approx 0$. Since $\alpha_M \geq \frac{1}{c}M^d$ and $\alpha_N \leq cN^d$, this gives us that $\frac{M^d}{c^2N^d} \leq \frac{\alpha_M}{\alpha_N} \approx 0$. Since $\frac{1}{c^2}$ is an element of $\bR$, this means $(\frac{M}{N})^d \approx 0$ and thus $\frac{M}{N}\approx 0$. 
\end{proof}

\subsection{Asymptotic cones for finitely generated groups}

In this section, we recall the \textbf{asymptotic cone} construction for finitely generated groups, using the approach pioneered by van den Dries and Wilkie \cite{LouAlex} via nonstandard analysis. (While asymptotic cones can be defined for any metric space, we will only need the construction for finitely generated groups equipped with their word metric.)

Throughout this subsection, we fix a finitely generated group $G$ with finite (symmetric) generating set $S$.  Given $N \in $ $^*\bN \backslash \bN$, we define the relation $\sim_N$ on $^*G$ by setting $x \sim_N y$ if and only if $\frac{d_S(x,y)}{N}$ is infinitesimal. The following properties of $\sim_N$ are routine to verify:

\begin{lem}\label{basicequiv}
Fix $M,N\in {}^*\bb N\setminus \bb N$ and $x,y,z\in {}^*G$.
\begin{enumerate}
\item $\sim_N$ is an equivalence relation on $^*G$.
\item Let $S'$ be another finite generating set for $G$, with corresponding word metric $d_{S'}$ and equivalence relation $\sim_N'$. Then $\sim_N$ and $\sim_N'$ are the same relation on $^*G$.
\item If $\frac{M}{N}$ is neither infinite nor infinitesimal, than $\sim_M$ and $\sim_N$ are the same equivalence relation. 
\item If $x \sim_N y$, then $\frac{|d_S(x,z) - d_S(y,z)|}{N}$ is infinitesimal. 
\item We have $x \sim_N y \iff y \in B(x,M)$ for some $M \in $ $^*\bN$ such that $\frac{M}{N} \approx 0$. 
\end{enumerate}
\end{lem}



For $x\in ^*G$ and $N\in {}^*\bb N\setminus \bb N$, we let $\operatorname{eq}_N(x)$ denote the equivalence class of $x$ with respect to $\sim_N$. Items (2) and (3) in Lemma \ref{basicequiv} imply that these equivalence classes are independent of the choice of generating set for $G$ and the choice of representative for the archimedean class of $N$.  Given $E\subseteq G$, we set $\operatorname{eq}_N(E)=\{\operatorname{eq}_N(x) \ : \ x\in E\}$ for the set of equivalence classes of elements of $E$.  

In the sequel, to simplify notation, we often overline a variable name to indicate that it is an equivalence class and similarly for sets of equivalence classes.  Consequently, for $\overline{x} \in \operatorname{eq}_N( ^*G)$, we have that $\overline{x}=eq^{-1}_N(\overline{x}) = \{g \in $ $^*G$ : $\operatorname{eq}_N(g) = \overline{x}\}$.  Implicit in this notation is that $x$ itself is a representative of the equivalence class $\overline{x}$, that is, $\operatorname{eq}_N(x)=\overline{x}$.  (One drawback of this notation is that we suppress mention of $N$; we thus only use this notation when $N$ has been fixed ahead of time.)

\begin{defn}
The \textbf{asymptotic cone} $C(G, g, N)$ of $G$ with base point $g \in $ $^*G$ and scaling factor $N \in$ $^*\bN \backslash \bN$ is defined as
$$
\left\{\operatorname{eq}_N(x) : x \in {}^*G \text{ and } \frac{d_S(x,g)}{N} \text{ is finite}\right\}.
$$
\end{defn}

By Lemma \ref{basicequiv}(4), the definition above is indeed independent of choice of representative.  The asymptotic cone can be equipped with the metric $d_N$ given by
$$d_N\left(\operatorname{eq}_N(x), \operatorname{eq}_N(y)\right) = \operatorname{st}\left(\frac{d_S(x,y)}{N}\right).$$ Note that this standard part is indeed defined as $\frac{d_S(x,y)}{N} \leq \frac{d_S(x,g)}{N} + \frac{d_S(y,g)}{N}$ and both terms on the right hand side of this inequality are finite by assumption. Furthermore, Lemma \ref{basicequiv}(4) implies that this definition is indeed independent of the choice of representatives for equivalence classes.  

An easy saturation argument shows that $C(G,g,N)$ is a complete metric space. Moreover, the asymptotic cone is also a geodesic metric space.  (See \cite{Ultrafilters} for proofs of both of these facts.)

We have obscured mention of the generating set $S$ in the above notation.  The reason for that is the fact that any two word metrics for $G$ are bi-Lipschitz equivalent remains true for the asymptotic cones with respect to these generating sets (using the same base point and scaling factor $N$), whence the construction is essentially independent of the choice of generator.  The map $$\operatorname{eq}_N(x)\mapsto \operatorname{eq}_N(gx):C(G,e,N)\to C(G,g,N)$$ is an isometry, whence the asymptotic cone (with respect to a given $N$) is also essentially independent of base point.  That being said, the homeomorphism type of the asymptotic cone can very much depend on the choice of $N$; see, for example, \cite{kramer}.  However, for the groups we will be interested in, the choice of $N$ will be irrelevant.

\begin{example}
All asymptotic cones of $\bb Z$ with respect to the generating set $\{1,-1\}$ are isometric to $\bb R$ with its usual distance.  Indeed, working for simplicity with the base point $0$, we see that the map $$\operatorname{eq}_N(x)\mapsto \operatorname{st}(\frac{x}{N}):C(\bb Z,0,N)\to \bb R$$ is an isometry.
\end{example}

  The result in the previous example is a special case of a more general result of Pansu \cite{PansuCones}, namely that when $G$ is virtually nilpotent, the isometry type of the asymptotic cones of $G$ (with respect to a given generating set) is independent of the choice of $N\in {}^*\bb N\setminus \bb N$.
 


The following theorem will be of great use to us going forward:

\begin{thm}
\label{propercone}
If $X$ is the metric space of an asymptotic cone defined via a finitely generated group $G$, then the following are equivalent:
\begin{enumerate}
    \item $X$ is proper (that is, closed balls are compact).
    \item $X$ is separable.
    \item $G$ is virtually nilpotent.
\end{enumerate}
In particular, some asymptotic cone of $G$ is proper (resp. separable) if and only if they all are.
\end{thm}

\begin{proof}[Proof sketch]
(1) implies (2) holds for any metric space; the converse was shown by Sisto in \cite{Sisto}.

(1) implies (3) was shown by Sapir \cite{Sapir}, building on work by Hrushovski \cite{hrushovski}.  Conversely, if $G$ is nilpotent, then some asymptotic cone of $G$ is proper by the work of van den Dries and Wilkie \cite{LouAlex}; that all such cones must be proper follows from the aforementioned result of Pansu \cite{PansuCones}.
\end{proof}

\section{A measure on the unit ball of asymptotic cones}

\subsection{Introducing the measure}

Throughout this section, $G$ is a finitely generated group with finite (symmetric) generating set $S$.  

For $N \in $ $^*\bN \backslash \bN$ and $g \in $ $^*G$, we let $\mu_{g,N}$ denote the Loeb measure on the hyperfinite set $B(g,N)$ and let $\cL_{g,N}$ denote the $\mu_{g,N}$-measurable subsets of $B(g,N)$. Let $\cA_{g,N}$ denote the set of subsets  $\overline{A}$ of $\operatorname{eq}_N\big(B(g,N)\big)$ such that $\overline{A} \in \cA_{g,N}$ if and only if $\operatorname{eq}_N^{-1}(\overline{A}) \cap B(g,N) \in \cL_{g,N}$. For $\overline{A} \in \cA_{g,N}$, set $\lambda_{g,N}(\overline{A}) = \mu_{g,N}\big(\operatorname{eq}_N^{-1}(\overline{A}) \cap B(g,N)\big)$. In other words, $\cA_{g,N}$ and $\lambda_{g,N}$ are the pushforwards of $\cL_{g,N}$ and $\mu_{g,N}$ by $\operatorname{eq}_N$.

Suitably adapting the example from Subsection 2.1, letting $G = \bZ$, $g = 0$, and $N$ be any element of $^*\bN\setminus \bb N$, $\lambda_{g,N}$ is the Lebesgue measure on $\lbrack -1,1 \rbrack$ and $\cal A_{g,N}$ is the collection of Lebesgue measurable subsets of $[-1,1]$. 

Note that, in the metric space $(C(G,g,N), d_N)$, $\operatorname{eq}_N\big(B(g,N)\big)$ is the closed ball of radius $1$ centered at the equivalence class $\overline{g} = \operatorname{eq}_N(g)$. In what follows, we will write $B_N(\overline{g},r)$ to refer to $\operatorname{eq}_N\big(B(g,rN)\big)$, which is the closed ball of radius $r \in \bR$ with center $\overline{g}$ in the asymptotic cone $C(G,g,N)$. Note that this is independent of the choice of representative of $\overline{g}$. We also consider the open ball of radius $r$ and center $\overline{g}$, denoted $U_N(\overline{g},r)$, and note that $U_N(\overline{g},r) = \bigcup_{n \in \bN} \operatorname{eq}_N\left(B\big(g,r(1-\frac{1}{n})N\big)\right)$.  Finally, we also consider the sphere $S_N(\overline{g},r) = B_N(\overline{g},r) \backslash U_N(\overline{g},r)$.  
 
 
 \begin{prop}\label{internalaremeasurable}
 For any $g \in$  $^*G$ and any $N \in $ $^*\bN\setminus \bb N$, if $A \subseteq $ $^*G$ is internal, then $\operatorname{eq}_N\big(A \cap B(g,N)\big)$ is $\lambda_{g,N}$-measurable.
 \end{prop}
 
 \begin{proof}
Set $E = A \cap B(g,N)$. It suffices to establish that $$\operatorname{eq}_N^{-1}\left(\operatorname{eq}_N(E)\right) = \bigcap_{n=1}^{\infty} \left\{x \in B(g,N) : x \in B\left(y,\frac{N}{n}\right)\text{ for some }y \in E\right\}.$$ 
 
 Towards this end, first assume $x \in \operatorname{eq}_N^{-1}\big(\operatorname{eq}_N(E)\big)$. Then there is $y \in E$ such that $x \in B(y,M)$ for some $M\in ^*\bb N$ with $\frac{M}{N} \approx 0$. Thus $x \in B(y,\frac{N}{n})$ for all $n \in \bN$, as desired.
 
  Conversely, assume that $x$ is in the intersection on the right side of the above display. Then for any $n \in \bN$, there is $a_n \in E$ such that $x \in B(a_n, \frac{N}{n})$. Note that $\{M \in $ $^*\bN : x \in B(a,M)$ for some $a \in E\}$ is internal since $E$ is internal, and thus has a minimum $M'$. By the above observation, $M' < \frac{N}{n}$ for all $n \in \bN$, meaning $\frac{M'}{N} \approx 0$. Taking $a \in E$ such that $x \in B(a,M')$, we see that $\operatorname{eq}_N(x)=\operatorname{eq}_N(a)$ and thus $x \in \operatorname{eq}_N^{-1}\big(\operatorname{eq}_N(E)\big)$.
\end{proof}

 \begin{cor}
 If $B_N(\overline{g},1)$ is separable (that is, if $G$ is virtually nilpotent), then the Borel subsets of $B_N(\overline{g},1)$ are $\lambda_{g,N}$-measurable.
 \end{cor} 
 
 \begin{proof}
First note that, since $B_N(\overline{g},1)$ is separable, the Borel subsets form the smallest $\sigma$-algebra containing all open balls, whence it suffices to show that open balls are $\lambda_{g,N}$-measurable. By Proposition \ref{internalaremeasurable}, the closed ball $B_N(\overline{h},r)$ is $\lambda_{g,N}$-measurable as $B_N(\overline{h},r)=\operatorname{eq}_N(B(h,rN)\cap B(g,N))$. Using the fact that $U_N(\overline{h},r) = \bigcup_{n \in \bN} B_N\big(\overline{h}, (1-\frac{1}{n})r\big)$, we see that all open balls are $\lambda_{g,N}$-measurable as well. 
\end{proof} 

From now on, we will restrict ourselves to the virtually nilpotent case and will only be concerned with the behavior of $\lambda_{g,N}$ on the Borel subsets of $B_N(\overline{g},1)$; we consequently restrict our attention to that a priori smaller $\sigma$-algebra. 

\subsection{The Lebesgue Density Property}

 One property that will come up frequently in our later proofs is whether or not $\lambda_{g,N}$ satisfies a version of the Lebesgue density theorem, which we formulate in the following way. 
 \begin{defn}
 Let $\mu$ be an outer measure on the metric space $X$.  We say that $(X,\mu)$ satisfies the \textbf{Lebesgue Density Theorem} (LDT) if, for all $\mu$-measurable $A\subseteq X$, we have that 
$$
\lim_{r \to 0^+} \frac{\mu\big(A \cap B(x,r)\big)}{\mu\big(B(x,r)\big)} = 1
$$
for $\mu$-almost all $x$ in $A$. 
 \end{defn}
 
 In order to put ourselves in the context of the previous definition, we define the outer measure $\lambda^0_{g,N}$ on $B_N(\overline{g},1)$ by $$\lambda^0_{g,N}(A) = \inf \{\lambda_{g,N}(B) : A \subseteq B, B \text{ a Borel set}\}.$$ This is an outer measure which agrees with $\lambda_{g,N}$ on Borel sets and whose measurable sets include the Borel sets. 
 
 
 The goal of this subsection will be to establish the following theorem:
 
 \begin{thm}\label{mainLDT}
 Let $G$ be a finitely generated virtually nilpotent group. Then for all infinite $N \in $ $^*\bN$ and all $g \in $ $^*G$, we have that $\lambda_{g,N}^0$ satisfies the Lebesgue Density Theorem. 
 \end{thm}
 
 

If we let $G = \bZ$ with the standard generating set, we have already discussed the fact that, for any $g \in {}^*\bZ$, $\lambda_{g,N}$ is the Lebesgue measure on  $B_N(\overline{g},1) \cong \lbrack -1, 1 \rbrack$.  Thus, the previous theorem generalizes the basic fact that the Lebesgue measure on $[-1,1]$ satisfies the conclusion of the Lebesgue Density Theorem.


Our proof of Theorem \ref{mainLDT} will use a more general result about the validity of the LDT for outer measures satsisfying a certain list of properties.

\begin{defn}
Suppose that $\mu$ is an outer measure on the metric space $X$.  We say that $\mu$:
\begin{enumerate}
    \item is \textbf{Borel regular} if all Borel sets are $\mu$-measurable and, for all $A \subseteq X$, we have that $\mu(A) = \mu(B)$ for some Borel set $B$; 
    \item is \textbf{open} $\sigma$\textbf{-finite} if $X$ can be covered by countably many open sets of finite measure; 
    \item has the \textbf{symmetric Vitali property (SVP)} if, for any $E \subseteq X$ with $\mu(E) < \infty$ and any collection $\cB$ of closed balls with centers in $E$ satisfying $$\inf\{r \in \bR^+: B(x,r) \in \cB\} = 0$$ for every $x \in E$, there is a countable, pairwise disjoint subcollection $\cB' = \{B_n : n \in \bN\} \subseteq \cB$ of balls such that $\mu(E \backslash \bigcup_{n \in \bN} B_n) = 0$; 
    \item has the \textbf{doubling property} if there exists a constant $C$ with $\mu\big(B(x,r)\big) \leq C\mu\big(B(x,\frac{r}{2})\big)$ for all $x \in X$ and $r \in \bR^+$. 
\end{enumerate}
\end{defn}

The following lemma appears as \cite[Remark 3.11]{SimonText}. 

\begin{lem}
\label{Vitali}
If $X$ is separable and $\mu$ is a Borel regular outer measure on X with the doubling property, then $\mu$ has the Symmetric Vitali Property. 
\end{lem}

The following is \cite[Theorem 3.16]{SimonText}:

\begin{thm}
\label{Lebesgue}
Let $\mu$ be an outer measure on the metric space $X$ that  is Borel regular, open $\sigma$-finite, and has the symmetric Vitali property.
Then $(X,\mu)$ satisfies the Lebesgue Density Theorem.
\end{thm}

We are now ready for the proof of the main result of this subsection:

\begin{proof}[Proof of Theorem \ref{mainLDT}]
The fact that $\lambda_{g,N}^0$ is Borel regular and open $\sigma$-finite is clear.  It remains to verify the symmetric Vitali property. Using Lemma \ref{Vitali} and the fact that $B_N(\overline{g},1)$ is separable (by the assumption that $G$ is virtually nilpotent), it suffices to establish that $\lambda_{g,N}^0$ has the doubling property. 

First recall that, by transfer, there is $c>0$ and $d\in \bb N$ such that $$\frac{N^d}{c} \leq |B(g,N)| \leq cN^d.$$ 

\

\noindent \textbf{Claim:}  For $\overline{x} \in B_N(\overline{g},1)$ and $r \in (0,2)$, we have:
$$
\lambda_{g,N}\left(B_N\big(\overline{x},\frac{r}{2}\big)\right) \geq \frac{\frac{1}{c}(\frac{r}{16}N)^d}{|B(g,N)|}.
$$

\

\noindent \textbf{Proof of claim:}  

There are four cases to consider:

\

\noindent Case 1:  Assume $B_N(\overline{x},\frac{r}{2}) \subseteq B_N(\overline{g},1)$. Then
$$
\lambda_{g,N}\left(B_N\big(\overline{x},\frac{r}{2}\big)\right) \geq \frac{\big|B(x,\frac{r}{4}N)\big|}{|B(g,N)|} \geq \frac{\frac{1}{c}(\frac{r}{4}N)^d}{|B(g,N)|} \geq \frac{\frac{1}{c}(\frac{r}{16}N)^d}{|B(g,N)|}.
$$

In the remaining three cases, we assume $B_N(\overline{x},\frac{r}{2})$ is not a subset of $B_N(\overline{g},1)$ and set $a := d_N(\overline{x},\overline{g}) > 1-\frac{r}{2}$. 

\noindent Case 2:  Assume $r \leq 1$. Then, using the fact that the asymptotic cone is a geodesic metric space, we can take a $\overline{z}$ such that $d_N(\overline{x},\overline{z}) = \frac{r}{4}$ and $d_N(\overline{z},\overline{g}) = a-\frac{r}{4} > 0$, whence $B_N(\overline{z},\frac{r}{4})\subseteq B_N(\overline{g},1)\cap B_N(\overline{x},\frac{r}{2})$. Thus, we can use the inequality from the first case, replacing $r$ with $\frac{r}{2}$, giving us
$$
\lambda_{g,N}\left(B_N\big(\overline{x},\frac{r}{2}\big)\right) \geq \lambda_{g,N}\left(B_N\big(\overline{z},\frac{r}{4}\big)\right) \geq \frac{\big|B(z,\frac{r}{8}N)\big|}{|B(g,N)|} \geq \frac{\frac{1}{c}(\frac{r}{16}N)^d}{|B(g,N)}.
$$

 In the last last two cases, we assume $r > 1$. 
 
 \noindent Case 3:  Assume $a \geq \frac{1}{2}$. Then we can take $\overline{z}$ with $d_N(\overline{x},\overline{z}) = \frac{1}{4}$ and $d_N(\overline{z},\overline{g}) = a-\frac{1}{4} > 0$, whence $B_N(\overline{z},\frac{1}{4})\subseteq B_N(\overline{x},\frac{r}{2})\cap B_N(\overline{g},1)$. Thus, we can use the inequality from the first case with $r = \frac{1}{2}$, giving us
 $$
 \lambda_{g,N}\bigg(B_N\Big(\overline{x},\frac{r}{2}\Big)\bigg) \geq \lambda_{g,N}\left(B_N\Big(\overline{z},\frac{1}{4}\Big)\right) \geq \frac{\big|B(z,\frac{1}{8}N)\big|}{|B(g,N)|} \geq \frac{\frac{1}{c}(\frac{1}{8}N)^d}{|B(g,N)|} \geq \frac{\frac{1}{c}(\frac{r}{16}N)^d}{|B(g,N)|}.
 $$
 
\noindent Case 4: $a < \frac{1}{2}$. Then there exists $\overline{z}$ with $d_N(\overline{x},\overline{z}) = \frac{1}{4}$, so that $d_N(\overline{z},\overline{g}) < \frac{3}{4}$. Thus $B_N(\overline{z},\frac{1}{4})$ is again contained in both of our balls, and so we can use the same inequality from case three. The claim is now proven.

Next note that, for any $r>0$, we have 

$$
\lambda_{g,N}\big(B_N(\overline{x},r)\big) \leq \frac{|B(x,2rN)|}{|B(g,N)|} \leq \frac{c(2rN)^d}{|B(g,N)|}.
$$
Combining this observation with the claim gives, for all $r \in (0,2)$, 

$$
\lambda_{g,N}\big(B_N(\overline{x},r)\big) \leq \frac{c(2rN)^d}{|B(g,N)|} = \frac{c^2(32)^d\frac{1}{c}(\frac{r}{16}N)^d}{|B(g,N)|} \leq c^2(32)^d\lambda_{g,N}\left(B_N\big(\overline{x},\frac{r}{2}\big)\right),
$$
and so the doubling property holds with $C = c^2(32)^d$ (as $\lambda_{g,N}^0$ agrees with $\lambda_{g,N}$ on Borel sets).
\end{proof}

\begin{remark}
Besides showing that there is a well-defined isometry type of asymptotic cones for a given virtually nilpotent group, Pansu \cite{PansuCones} also showed that this asymptotic cone can be equipped with a natural Lie group structure.  (See \cite{SphereGrowth} for a nice description of this construction.)  In particular, these asymptotic cones possess an (essentially) unique Haar measure; we leave it for future work to compare the restriction of this Haar measure to the unit ball of the asymptotic cone and the measure on the unit ball we have been considering in this paper.  (We thank Alessandro Sisto for pointing us to this fact.)
\end{remark}


\subsection{Consequences of the small sphere property}

While the assumption of polynomial growth will prove useful going forward (as has already been evidenced by the validity of the LDT for the measure we have considered on closed balls in asymptotic cones established in the previous subsection), it does not seem to be quite powerful enough to prove all of our results below. Indeed, we will want spheres in asymptotic cones to behave much as they do in euclidean space, namely that a sphere of radius $r$ of dimension (that is, degree of polynomial growth) $d$ has size on the order of $r^{d-1}$.  We formalize this desire into the following definition:

\begin{defn}
A finitely generated group $G$ with symmetric generating set $S$ has the \textbf{Small Sphere (SS) Property} if there exists some constant $c'$ such that, for any $n \in \bN$ and $g \in G$, we have that $|S(g,n)| \leq c'n^{d-1}$.
\end{defn}

As defined, the SS property is a property that may or may not hold for a given group \emph{with respect to a given generating set}; it is unclear to us if this property is independent of generating set. 

The SS property has been shown to hold in nilpotent groups of nilpotency class 2 (such as the discrete Heisenberg group), regardless of the generating set, by Stoll in \cite{2StepStoll}. Indeed, Stoll was able to show that this growth rate of spheres is sharp in this case, that is, there exist nilpotent groups of nilpotecy class 2 where the size of spheres of radius $n$ is precisely on the order of $n^{d-1}$. 

In \cite{SphereGrowth}, Breuillard and Le Donne made progress on finding a growth rate for spheres in arbitrary virtually nilpotent groups. They proved that for a nilpotent group of nilpotency class $r$, the size of the spheres of radius $n$ is bounded by a growth function on the order of $n^{d-\beta_r}$, where $\beta_r = \frac{2}{3r}$. However, they mention that this growth rate is unlikely to be sharp. Indeed, it has been conjectured that spheres in all  virtually nilpotent groups satisfy the SS property. Breuillard and Le Donne refer to this as a ``folklore conjecture'' and no further progress on this conjecture seems to have been made at this time. 

Nevertheless, we will concern ourselves primarily with groups with the SS property (with respect to some generating set), as this will allow us to prove many of the nice properties of SIM subsets of $\bZ$ in our more general context. One of the main consequences this assumption gives us is the following:

\begin{lem}
\label{NiceBalls}
Suppose that $G$ is a finitely generated virtually nilpotent group with the SS property. Fix  $N \in$ $^*\bN \backslash \bN$ and $g \in $ $^*G$. Then $\lambda_{g,N}\big(S_N(\overline{x},r) \cap B_N(\overline{g},1)\big) = 0$ for all $r \in \bR^+$ and all $x \in $ $^*G$. 
\end{lem}

\begin{proof}
For $\delta > 0$, let $B_{\delta} = \bigcup_{n \in I_{\delta}} S(x,n)$, where $I_{\delta} \subseteq $ $^*\bN$ is the infinite, hyperfinite interval $\big\lbrack \lfloor (1-\delta)rN \rfloor +1 , \lfloor (1+\delta)rN \rfloor \big\rbrack$. Note that $B_{\delta}$ is internal and $S_N(\overline{x},r) \subseteq B_{\delta}$, whence $\lambda_{g,N}\big(S_N(\overline{x},r) \cap B_N(\overline{g},1)\big) \leq \frac{|B_{\delta} \cap B(g,N)|}{|B(x,N)|} \leq \frac{|B_{\delta}|}{|B(x,N)|}$. Assume that $c$ is the constant witnessing the polynomial growth of $G$ with degree $d$, and $c'$ is the corresponding constant for the SS property. Then we have that $|I_{\delta}| \leq 2\delta rN$, and $|S(x,M)| \leq c'\big((1+\delta)rN\big)^{d-1}$ for all $M \in I_{\delta}$. Thus, $|B_{\delta}| \leq (2\delta rN)\big(c'\left((1+\delta)rN\right)^{d-1}\big)$ and $|B(x,N)| \geq \frac{1}{c}N^d$. Therefore,
$$
\lambda_{g,N}\big(S_N(\overline{x},r) \cap B_N(\overline{g},1)\big) \leq \frac{2c'(1+\delta)^{d-1}\delta r^dN^d}{\frac{1}{c}N^d}  = 2cc'(1+\delta)^{d-1}\delta r^d
$$
for all $\delta > 0$. Letting $\delta$ tend to $0$, we see that $\lambda_{g,N}\big(S_N(\overline{x},r) \cap B_N(\overline{g},1)\big) = 0$, as desired. 
\end{proof}

One nice consequence of Lemma \ref{NiceBalls} is that the measure of balls becomes exactly what we would expect them to be, as long as they are of reasonable size.
\begin{cor}
\label{NiceMeasures}
Suppose that $G$ is a finitely generated virtually nilpotent group with the SS property.  Then for any $N\in \bb {}^*N\setminus \bb N$, $r \in \bb R^{^+}$, and $g,h \in $ $^*G$, we have that  
$$
\lambda_{g,N}\Big(\operatorname{eq}_N\big(B(h,rN)\cap B(g,N)\big)\Big) = \operatorname{st}\left(\frac{\big|B(h,rN) \cap B(g,N)\big|}{|B(g,N)|}\right).
$$
\end{cor}

\begin{proof}
By Lemma \ref{NiceBalls}, we know that $\lambda_{g,N}\left(\operatorname{eq}_N\big(S(h,rN)\big)\cap B_N(\overline{g},1)\right) = 0$ and $\lambda_{g,N}\left(\operatorname{eq}_N\big(S(g,N)\big)\right) = 0$ (as $\operatorname{eq}_N\big(S(h,rN)\big) \subseteq S_N(\overline{h},r)$ and $\operatorname{eq}_N\big(S(g,N)\big) \subseteq S_N(\overline{g},1)$). Therefore, the left hand side of the display appearing in the corollary equals
\small
$$
 \lambda_{g,N}\Big(\operatorname{eq}_N\big(B(h,rN)\cap B(g,N)\big) \backslash \big(\operatorname{eq}_N\big(S(h,rN)\big)\cup \operatorname{eq}_N\big(S(g,N)\big)\big)\Big).
$$
\normalsize

 Notice that $\operatorname{eq}_N^{-1}\Big(\operatorname{eq}_N\big(B(h,rN)\cap B(g,N)\big) \backslash \left(\operatorname{eq}_N\big(S(h,rN)\big) \cup \operatorname{eq}_N\big(S(g,N)\big)\right)\Big)$ is a subset of both $B(h,rN)$ and $B(g,N)$, since the only equivalence classes left would be those that were some appreciable distance from $S(h,rN)$ and $S(g,N)$. Consequently, we have 
$$
\lambda_{g,N}\Big(\operatorname{eq}_N\big(B(h,rN)\cap B(g,N)\big)\Big) \leq \mu_{g,N}\big(B(h,rN) \cap B(g,N)\big).
$$ 

On the other hand, $B(h,rN) \cap B(g,N)$ is a subset of $\operatorname{eq}_N^{-1}\left(\operatorname{eq}_N\big(B(h,rN) \cap B(g,N)\big)\right)$, meaning $\lambda_{g,N}\left(\operatorname{eq}_N\big(B(h,rN) \cap B(g,N)\big)\right) \geq \mu_{g,N}\big(B(h,rN) \cap B(g,N)\big)$. Combining these bounds, we see that the left-hand side of the display in the corollary equals
$$ \mu_N\big(B(h,rN)\cap B(g,N)\big) = \operatorname{st}\left(\frac{\big|B(h,rN) \cap B(g,N)\big|}{|B(g,N)|}\right),
$$ 
as desired.
\end{proof}

Using the same techniques, we prove the following lemmas hold in virtually nilpotent groups with the SS property. The first lemma tells us that we can restrict internal sets to our given ball before or after moving to the equivalence classes without changing the measure.

\begin{lem}
\label{beforeandafter}
Suppose that $G$ is a finitely generated virtually nilpotent group with the SS property.  Then for any $g \in $ $^*G$, any $N \in $ $^*\bN \backslash \bN$, and any internal $A \subseteq $ $^*G$, we have
$$
\lambda_{g,N}\big(\operatorname{eq}_N(A) \cap \operatorname{eq}_N(B(g,N))\big) = \lambda_{g,N} \left(\operatorname{eq}_N\big(A \cap B(g,N)\big)\right).
$$
\end{lem}

\begin{proof}
Since $\operatorname{eq}_N\big(A \cap B(g,N)\big) \subseteq \operatorname{eq}_N(A) \cap \operatorname{eq}_N(B(g,N))$, the left-hand side is at least the right-hand side. To show equality, note that 
$$
\big(\operatorname{eq}_N(A) \cap \operatorname{eq}_N(B(g,N))\big) \backslash \operatorname{eq}_N\big(A \cap B(g,N)\big) \subseteq S_N(\overline{g},1),
$$  
as any equivalence class in $\operatorname{eq}_N(A)$ and $\operatorname{eq}_N(B(g,N))$ that is not represented by an element of $A \cap B(g,N)$ must contain elements outside of $B(g,N)$. Thus, by Lemma \ref{NiceBalls}, 
$$
\lambda_{g,N}\big(\operatorname{eq}_N(A) \cap \operatorname{eq}_N(B(g,N))\big) = \lambda_{g,N}\left(\big(\operatorname{eq}_N(A) \cap \operatorname{eq}_N(B(g,N))\big) \backslash S_N(\overline{g},1)\right) = 
$$
$$
= \lambda_{g,N}\left(\operatorname{eq}_N\big(A \cap B(g,N)\big)\right). \text{ \qedhere}
$$
\end{proof}

This next lemma gives us an easy way to convert the measure from one ball to a superball or subball. 

\begin{lem}
\label{Subballs}
Suppose that $G$ is a finitely generated virtually nilpotent group with the SS property.  Then for any $g, g' \in $ $^*G$, any $N \in $ $^*\bN \backslash \bN$, any $r \in (0,2)$ such that $B(g',rN) \subseteq B(g,N)$, and any internal $A \subseteq $ $^*G$, we have
$$
\lambda_{g,N}\Big(\operatorname{eq}_N\big(A \cap B(g',rN) \big)\Big) = \operatorname{st}\left(\frac{|B(g',rN)|}{|B(g,N)|} \right)\lambda_{g',rN}\Big(\operatorname{eq}_N\big(A \cap B(g',rN)\big)\Big).
$$
\end{lem}

\begin{proof}
By definition, whenever $\overline{B} \subseteq \operatorname{eq}_N(^*G) = \operatorname{eq}_{rN}(^*G)$ is $\lambda_{g,N}$ or $\lambda_{g',rN}$ measurable, respectively, we have $\lambda_{g,N}(\overline{B}) = \mu_{g,N}\big(\operatorname{eq}_N^{-1}(\overline{B}) \cap B(g,N)\big)$ and $\lambda_{g',rN}(\overline{B}) = \mu_{g',rN}\big(\operatorname{eq}_N^{-1}(\overline{B}) \cap B(g',rN)\big)$. Next observe that 
$$
\Big(\operatorname{eq}_N^{-1}\big(\operatorname{eq}_N\big(B(g',rN)\big)\big) \cap B(g,N)\Big) \Big\backslash \Big(\operatorname{eq}_N^{-1}\big(\operatorname{eq}_N\big(B(g',rN)\big)\big) \cap B(g',rN) \Big) 
$$
 is contained in $\operatorname{eq}_N^{-1}(S_N(\overline{g'},r))$, which has $\mu_{g,N}$-measure $0$ by Lemma \ref{NiceBalls}, and whence the restriction to $B(g',rN)$ does not affect the measure. The lemma follows from noting that the following four quantities are equal, where the infima in the second and third bullets are over internal subsets $C$ of $^*G$:
\begin{itemize}
\item $\mu_{g,N}\Big(\operatorname{eq}_N^{-1}\Big(\operatorname{eq}_N\big(A \cap B(g',rN)\big)\Big) \cap B(g,N)\Big)$ 
\item $\inf\left\{\operatorname{st}\left(\frac{|C|}{|B(g,N)|}\right) : \operatorname{eq}_N^{-1}\Big(\operatorname{eq}_N\big(A \cap B(g',rN)\big)\Big) \cap B(g,N) \subseteq C\right\}$
\item $\inf\left\{\operatorname{st}\left(\frac{|C|}{|B(g',rN)|}\cdot \frac{|B(g',rN)|}{|B(g,N)|}\right) : \operatorname{eq}_N^{-1}\Big(\operatorname{eq}_N\big(A \cap B(g',rN)\big)\Big) \cap B(g',rN) \subseteq C\right\}$ 
\item $\operatorname{st}\left(\frac{|B(g',rN)|}{|B(g,N)|}\right) \mu_{g',rN} \Big(\operatorname{eq}_N^{-1}\Big(\operatorname{eq}_N\big(A \cap B(g',rN)\big)\Big) \cap B(g',rN)\Big).$ \qedhere \end{itemize}
\end{proof}


\section{The BM property for internal sets}

\subsection{The Ball Measure Property}
Throughout this section, we fix a finitely generated group $G$ with symmetric generating set $S$.

\begin{defn}
Given a ball $I = B(g,N)$, with $g \in $ $^*G$ and $N \in $ $^*\bN \backslash \bN$, and internal $A \subseteq$ $^*G$, we define the following two quantities:
\begin{itemize}
    \item $\lambda_I(A) = \lambda_{g,N}\big(\operatorname{eq}_N(A \cap I)\big)$ .
    \item $g_A(I) = \max \left\{\frac{|B(x,M)|}{|I|} : B(x,M) \subseteq I,\text{ } B(x,M) \cap A = \emptyset\right\}$.
\end{itemize}
\end{defn} 
Using the notation from the previous definition, note that $\lambda_I(A)$ can equivalently be expressed as $\lambda_{e,N}\left(\operatorname{eq}_N\big(g^{-1}A \cap B(e,N)\big)\right)$.  We will sometimes refer to $\lambda_I(A)$ as the "measure" of $A$ in $I$. The quantity $g_A(I)$ is a measure of the size of gaps of $A$ in $I$ and notice that the maximum appearing in the definition indeed exists as $\big\{M \in $ $^*\bN : B(x,M) \subseteq I,$ $ B(x,M) \cap A = \emptyset \text{ for some } x \in I\big\}$ is an internal subset of $\lbrack 1, 2N \rbrack$.  Just like the corresponding properties in the original definition of the IM property, $\lambda_I(A)$ is an external property of $A$ whereas $g_A(I)$ is an internal property. As in the case of the integers, without any further assumptions, there is always one relationship between these two quantities, assuming the group is virtually nilpotent and satisfies the SS property:

\begin{lem}
\label{ReverseBM}
Suppose that $G$ is a virtually nilpotent group satisfying the SS property.  Let $I = B(g,N)$ be an infinite, hyperfinite ball and $A$ an internal subset of $I$. For any $\varepsilon > 0$, if $\lambda_I(A) \geq 1 - \varepsilon$, then $g_A(I) \leq \varepsilon$. 
\end{lem}

\begin{proof}
Assume that $g_A(I) > \varepsilon$. Then there is some infinite hyperfinite ball $J \subseteq I$ such that $\operatorname{st}(\frac{|J|}{|I|}) > \varepsilon$ and $A \cap J = \emptyset$. Let $\partial J$ be the sphere that makes up the boundary of $J$. Observe that $\operatorname{eq}_N(A) \cap \big(\operatorname{eq}_N(J) \backslash \operatorname{eq}_N(\partial J)\big) = \emptyset$. Thus, using Lemma \ref{NiceBalls} and Corollary \ref{NiceMeasures}, we have that the following quantities coincide:
\begin{itemize}
    \item $\lambda_I(A)$
    \item $\lambda_{g,N}\big(\operatorname{eq}_N(A)\big)$
    \item $ 1 - \lambda_{g,N}\big(\operatorname{eq}_N(J) \backslash \operatorname{eq}_N(\partial J)\big)$
    \item $1- \lambda_{g,N}\big(\operatorname{eq}_N(J)\big) + \lambda_{g,N}\big(\operatorname{eq}_N(\partial J)\big)$
    \item $1 - \lambda_{g,N}\big(\operatorname{eq}_N(J)\big)$
    \item $ 1 - \operatorname{st}\left(\frac{|J|}{|I|}\right)$.
\end{itemize}
Since the last quantity is strictly less than $1 - \varepsilon$, the lemma is proven.
\end{proof}

Notice that even this proof of such a basic relationship (seemingly) requires the use of the SS property.  The definition of the BM property seeks to find a connection between these quantities in the other direction:

\begin{defn}
Given an internal set $A \subseteq $ $^*G$ and infinite hyperfinite ball $I$, we say that $A$ has the \textbf{ball measure property} (or \textbf{BM property} for short) on $I$ if: for every $\varepsilon > 0$, there is $\delta > 0$ such that, for all infinite hyperfinite balls $J \subseteq I$ with $g_A(J) \leq \delta$, we have $\lambda_J(A) \geq 1 - \varepsilon$. 
\end{defn}

For simplicity, we often just say $A$ is BM on $I$. Note that if $G = \bZ$ with generating set $S = \{-1,1\}$, then the notion of the BM property reduces to the definition of the IM property on intervals of odd length which, since a single element has no effect on either $g_A(I)$ or $\lambda_I(A)$, can be easily extended to include all intervals. Thus, the previous definition truly is a generalization of the IM property. 

While we have defined this concept for all finitely generated groups, given that the basic fact Lemma \ref{ReverseBM} used the assumptions of being virtually nilpotent and having the SS property, it seems that to establish a nice theory of sets with the BM property, these assumptions will need to be present throughout. Consequently, \textbf{we henceforth assume that $G$ is of polynomial growth and has the SS property}.

\subsection{Properties of BM sets}

	Suppose that $A$ has the BM property on $I$.  Let $\delta(A,I,\varepsilon)$ be the supremum of the $\delta$'s that witness $A$ being BM on $I$ for the given $\varepsilon$. As in the IM case, we want to avoid the situation where $A$ is BM on $I$ simply because there is some $\delta > 0$ with $g_A(J) > \delta$ for all infinite balls $J \subseteq I$. In this spirit, we establish the following proposition:.
	
\begin{prop}
Let $A \subseteq $ $^*G$ be an internal set and let $I$ be an infinite hyperfinite ball such that $A$ is BM on $I$. Then the following statements are equivalent:
\begin{enumerate}
\item[(i)] For every $\delta > 0$, there is an infinite subball $J_{\delta}$ of $I$ such that $g_A(J_{\delta}) < \delta$.
\item[(ii)] There is an infinite subball $J$ of $I$ such that $\lambda_J(A) > 0$. 
\end{enumerate}
\end{prop}

\begin{proof}
First assume that (i) holds. Then there is a subball $J$ of $I$ such that $g_A(J) < \delta(A,I,\frac{1}{2})$. By definition, this means that $\lambda_J(A) \geq 1-\frac{1}{2} = \frac{1}{2} > 0$. 

For the other direction, assume that $J = B(g,N)$ is an infinite subball of $I$ such that $\lambda_J(A) > 0$, meaning $\lambda_{g,N}\big(\operatorname{eq}_N(A \cap J)\big) > 0$. By the LDT and Corollary \ref{NiceMeasures}, for any $\delta > 0$, we can find a subball $J' = B(g',N') \subseteq J$ such that 
$$
1-\delta < \frac{\lambda_{g,N}\big(\operatorname{eq}_N(A \cap J) \cap \operatorname{eq}_N(J')\big)}{\lambda_{g,N}(\operatorname{eq}_N(J'))} = \operatorname{st}\left(\frac{|J|}{|J'|}\right) \lambda_{g,N}\big(\operatorname{eq}_N(A \cap J) \cap \operatorname{eq}_N(J')\big). 
$$
By Lemma \ref{beforeandafter} and Lemma \ref{Subballs}, we have that the right-hand side of the previous display equals
$$
 \operatorname{st}\left(\frac{|J|}{|J'|}\right) \lambda_{g,N}\big(\operatorname{eq}_N(A \cap J')\big) = \lambda_{g',N'}\big(\operatorname{eq}_N(A \cap J')\big).
$$
 Thus, by Lemma \ref{ReverseBM}, we see that $g_A(J') < \delta$. Letting $J_{\delta} = J'$, we have that (i) holds.  
\end{proof}

Using the same terminology as in \cite{GoldbringBook}, we say that $A$ has the \textbf{enhanced BM property on $I$} if $A$ is BM on $I$ and $\lambda_I(A) > 0$. 

We can now establish the ``internal partition regularity'' of the BM property:

\begin{thm}
\label{PRBM}
Let $A$ be an internal subset of $^*G$ with the enhanced BM property on an infinite hyperfinite ball $I$. Assume that $A \cap I = B_1 \cup \cdots \cup B_n$ with each $B_i$ internal. Then there is an $i\in \{1,\ldots,n\}$ and an infinite subball $J \subseteq I$ such that $B_i$ has the enhanced BM property on $J$.
\end{thm}

\begin{proof}
We prove this theorem by induction, the base case $n=1$ being trivial. Now assume that $n>1$ and the result is true for partitions of size $n-1$. Let $A\cap I=B_1\cup \cdots\cup B_n$ be an internal partition. If there is an $i\in \{1,\ldots,n\}$ and an infinite subball $J \subseteq I$ such that $B_i \cap J = \emptyset$ and $\lambda_J(A) > 0$, then we can use the assumption on $J$ and the $n-1$ remaining $B_j$'s to obtain $j$ and $J' \subseteq J \subseteq I$ that satisfies the conclusion of the theorem. Thus, we may assume going forward that this is not the case, that is, whenever $\lambda_J(A) > 0$, we have that $B_i \cap J \neq \emptyset$ for all $i=1,\ldots,n$. We will show that this implies that each $B_i$ is BM on $I$. Since $\lambda_I(B_i)$ must be positive for some $i\in \{1,\ldots,n\}$, this $B_i$ wil be as desired. 

Fix $i$ and set $B = B_i$. Fix $\epsilon>0$ and let $J = B(g,N)$ be an infinite subball of $I$ with $g_B(J) \leq \delta(A,I,\varepsilon)$. (If no such $J$ exists then $B$ has the BM property on $I$ and we are done.)  Since $B \subseteq A$, we know that $g_A(J) \leq g_B(J) \leq \delta(A,I,\varepsilon)$. By definition, this gives us that $\lambda_J(A) \geq 1 - \varepsilon$. Now take $B_N(\overline{x},r) \subseteq B_N(\overline{g},1) \backslash \operatorname{eq}_N(B)$, so that $B(x,rN) \cap B = \emptyset$. Our assumption then gives us that $\lambda_{B(x,rN)}(A) = 0$, and so $\lambda_J\big(A \cap B(x,rN)\big) = 0$. Thus, $\lambda_J(B) = \lambda_J(A) \geq 1 - \varepsilon$, meaning $B$ has the $BM$ property on $I$. 
\end{proof}

We next proceed to establish a result about ``difference sets'' for which we will need the following notion:

\begin{defn}
Let $A_1,...,A_n$ be internal subsets of $^*G$ and let $I_1,...,I_n$ be hyperfinite balls. 
\begin{itemize}
    \item A $\delta$\textbf{-configuration} with respect to these sets is a sequence of subballs $J_1,...,J_n$ of $I_1,...,I_n$ respectively such that all $J_i$ have the same radius and such that $g_{A_i}(J_i) \leq \delta$ for all $i=1,\ldots,n$. We will refer to the shared radius of $J_1,...,J_n$ as the \textbf{radius of the configuration}.
    \item If $J_1,...,J_n$ is a $\delta$-configuration with respect to $A_1,...,A_n$ and $I_1,...,I_n$, then a $\delta$\textbf{-subconfiguration} of $J_1,\ldots,J_n$ is a $\delta$-configuration with respect to $A_1,...,A_n$ and $J_1,...,J_n$.
    \item A \textbf{strong $\delta$-subconfiguration} is a $\delta$-subconfiguration $K_1,...,K_n$ of $J_1,...,J_n$ for which there exists a $c \in $ $^*G$ such that, writing $J_i = B(a_i, N)$, we have that $K_i = B(a_ic,M)$.
\end{itemize}
\end{defn}

\begin{thm}
\label{DeltaConfigs}
Let $A_1,...,A_n$ be internal subsets of $^*G$ which are BM on the balls $I_1,...,I_n$ respectively. Take $\varepsilon > 0$ such that $\varepsilon < \frac{1}{n+1}$, and $\delta > 0$ such that $\delta < \min_{i=1,...,n} \delta(A_i,I_i,\varepsilon)$. Then there is a $w \in \bN$ such that any $\delta$-configuration has a strong $\delta$-subconfiguration of radius at most $w$. 
\end{thm}

\begin{proof}
Let $A_1,...,A_n$, $I_1,...,I_n$, $\varepsilon$, and $\delta$ be as above. Take $J_i = B(a_i,R)$ to be a $\delta$-configuration, with $R \in $ $^*\bN \backslash \bN$. Then, by our choice of $\delta$, we have that $\lambda_{J_i}(A_i) \geq 1-\frac{1}{n+1}$, meaning $\lambda_{a_i^{-1}J_i}(a_i^{-1}A_i) \geq 1- \frac{1}{n+1}$. Note that $a_i^{-1}J_i = B(e,R)$. Therefore, each $a_i^{-1}A_i$ has measure at least $1-\frac{1}{n+1}$ on $B(e,R)$, whence $$\lambda_{e,R}\left(\bigcap_{i=1}^n \operatorname{eq}_R\big(a_i^{-1}A_i \cap B(e,R)\big)\right) > 0.$$ By the Lebesgue density theorem, $\bigcap_{i=1}^n \operatorname{eq}_R\big(a_i^{-1}A_i \cap B(e,R)\big)$ has a point of density $\overline{b}$. Thus, there must be an $r \in (0,1)$ such that 

$$
\dfrac{\lambda_{e,R}\left(\Big(\bigcap_{i=1}^n \operatorname{eq}_R\big(a_i^{-1}A_i \cap B(e,R)\big)\Big) \cap B_R(\overline{b},r)\right)}{\lambda_{e,R}\big(B_R(\overline{b},r)\big)} \geq 1-\delta.
$$

Let $K_i = B\big(a_ib, \lfloor rR \rfloor \big)$.  By the previous display, we have that $g_{A_i}(K_i) \leq \delta$ for all $i=1,\ldots,n$. Thus the $K_i$ form a strong $\delta$-subconfiguration with radius $rR < R$. 

The import of the previous discussion is that every $\delta$-configuration with infinite radius has a \emph{proper} strong $\delta$-subconfiguration. Let $\cC$ denote the (internal) set of $\delta$-configurations of $A_1,...,A_n,I_1,...,I_n$ and let $f : \cC \to $ $^*\bN$ be the internal function defined by sending any $\delta$-configuration $J_1,...,J_n$ to the minimal radius of a strong $\delta$-subconfiguration of $J_1,...J_n$ (which exists by transfer). We have thus shown that $f(\cC) \subseteq \bN$ whence, by overflow, $f(\cC) \subseteq \lbrack 0, w \rbrack$ for some $w \in \bN$. Thus, this $w$ satisfies the conditions of the theorem. 
\end{proof}
 
The previous lemma has consequences for difference sets.  To state these consequences, we establish some terminology and remind the reader about syndetic sets in arbitrary groups.

\begin{defn}
For any $A \subseteq $ $^*G$, we define
$$
D(A) = \{g \in G : g = a^{-1}b \text{ for infinitely many pairs } a,b \in A\}.
$$
\end{defn}

For any group $G$, recall that a subset $A$ of $G$ is called \textbf{syndetic} if there is a finite subset $F \subseteq G$ such that $AF = G$. Equivalently, $A$ is syndetic if there exists $n \in \bN$ such that $B(g,n) \cap A \neq \emptyset$ for all $g \in G$.  Note that if $H \leq G$ is a subgroup, then $H$ is syndetic if and only if $H$ has finite index in $G$.

We can now state the following lemma, which we will use in the next section.

\begin{lem}
\label{BMSDA}
Suppose that $A$ has the enhanced BM property on some infinite, hyperfinite ball $I$. Then there is some $r \in \bN$ such that, for every $g \in G$, there is $z \in B(e,r)$ with $x^{-1}y \in \{h^{-1}ghz : h \in $ $^*G\}$ for infinitely many pairs $x,y \in A$. 
\end{lem}

\begin{proof}
Let $A_1 = A_2 = A$ and $I_1 = I_2 = I = B(v,N)$. By Theorem \ref{DeltaConfigs}, there is $w \in \bN$ such that every $\delta$-configuration $B(a_i,R)$ has some $c \in $ $^*G$ with the sequence $B(a_ic,w)$ containing a strong $\delta$-subconfiguration of $B(a_i,R)$. In particular, this implies that $B(a_ic,w) \cap A \neq \emptyset$. 

Fix $g \in G$. Since $\lambda_{v,N}\big(\operatorname{eq}_N(A \cap I)\big) > 0$, we can find countably many distinct points of density of $\operatorname{eq}_N(A \cap I)$, thereby giving us pairwise disjoint balls $B_N(\overline{a_{i,\varepsilon}},b_{\varepsilon}) \subseteq \operatorname{eq}_N(A \cap I)$ with $\lambda_{a_{i,\varepsilon},b_{\varepsilon}N}\big(\operatorname{eq}_N(A \cap I) \cap B_N(\overline{a_{i,\varepsilon}},b_{\varepsilon})\big) \geq 1-\varepsilon$ for all $i \in \bN$. Thus $g_A\big(B(a_{i,\varepsilon},b_{\varepsilon}N)\big) \leq \varepsilon$. Since the gap function is internal and we can find such balls for any $\varepsilon > 0$, there must be countably many pairwise disjoint balls $J_i = B(a_i,R)$ with $g_A(J_i) \approx 0$ for all $i\in \bb N$. Since $g \in G$, it is also true that $g_A\big(B(a_ig,R)\big)\approx 0$, as the elements of $B(a_ig,R)$ that are not in $J_i$ are all contained in one of $S(a_i,R+1),...,S(a_i,R + |g|)$, all of which have infinitesimal size when compared to $J_i$ by the SS property. Consequently, the $J_i, B(a_ig,R)$ forms a $\delta$-configuration of $A_1,A_2,I_1,I_2$ for any $\delta$ and any $i$. Thus, by our choice of $w$, for every $i \in \bN$, there is a $c_i \in $ $^*G$ such that
$$
A \cap B(a_ic_i,w) \neq \emptyset \text{ and } A \cap B(a_igc_i,w) \neq \emptyset.
$$

If $x_i \in A \cap B(a_ic_i,w)$ and $y_i  \in A \cap B(a_igc_i,w)$, then $x_i = a_ic_ix$ and $y_i = a_igc_iy$ with $|x|,|y| \leq w$. Since

$$
x_i^{-1}y_i = x^{-1}c_i^{-1}a_i^{-1}a_igc_iy = x^{-1}c_i^{-1}gc_ix(x^{-1}y),
$$

we see that $x_i^{-1}y_i\in B(h^{-1}gh,2w)$ for some $h \in $ $^*G$, meaning there must be a $z \in B(e,2w)$ such that infinitely many $i$ satisfy $x_i^{-1}y_i \in \{h^{-1}ghz: h \in $ $ ^*G\}$. 
\end{proof}

As a corollary, we get the following results.

\begin{cor}
\label{BMD(A)}
If $A$ has the enhanced BM property on some infinite, hyperfinite ball $I$, then there is some $r \in \bN$ such that $B(g,r) \cap D(A) \neq \emptyset$ for any $g \in Z(G)$. 
\end{cor}

\begin{proof}
Let $r$ be as in Lemma \ref{BMSDA}. Fix $g \in Z(G)$. By Lemma \ref{BMSDA}, there is $z \in B(e,r)$ such that there are infinitely many $x,y \in A$ with $x^{-1}y \in \{h^{-1}ghz : h \in $ $^*G\} = \{gz\}$. In other words, $gz \in D(A)$. Since $z \in B(e,r)$, $gz \in B(g,r)$, finishing the proof. 
\end{proof}

Thus, we see that sets with the enhanced BM property are, in some sense, "syndetic on the center of $G$".

\begin{cor}
If $A$ has the enhanced BM property on some infinite, hyperfinite ball $I$ and $Z(G)$ has finite index in $G$, then $D(A)$ is syndetic in $G$. 
\end{cor}

\section{The SBM property for standard sets}

\subsection{SBM Sets and Their Properties}

We now move to the setting of standard sets, \textbf{maintaining our standing assumptions that $G$ is a finitely generated virtually nilpotent group with the SS property}. 

\begin{defn}
We say that a set $A \subseteq G$ has the \textbf{standard BM property} (or \textbf{SBM property} for short) if:
\begin{enumerate}
    \item $^*A$ is BM on every infinite, hyperfinite subball of $^*G$, and 
    \item $^*A$ has the enhanced BM property on some infinite hyperfinite subball of $^*G$.
\end{enumerate}
\end{defn}

Note that for $G = \bZ$ with the standard generating set $S = \{1,-1\}$, we recover the definition of SIM set. 

Since the SBM property is a property of standard sets, it is possible to give an equivalent definition using only standard notions (as Leth did for subsets of the integers with the SIM property in \cite{LethSIM}).  We merely establish the notation for stating this equivalent reformulation, leaving the proof of the theorem to the reader.  Since this standard reformulation is not terribly enlightening, the reader choosing to skip the verification will lose nothing in what follows.

Following the notation in \cite{GoldbringBook}, for $A \subseteq G$ and $0 < \delta < \varepsilon < 1$, let the function $F_{\delta,\varepsilon,A} : \bN \to \bN \cup \{\infty\}$ be defined as follows. First, if $g_A(I) > \delta$ for every subball of $G$ with radius $\geq n$, set $F_{\delta,\varepsilon,A}(n) = 0$. Otherwise, let $F_{\delta,\varepsilon,A}(n)$ be the minimum $k$ such that there is a subball $I$ of $G$ with radius $\geq n$ and $g_A(I) \leq \delta$ for which there are disjoint subballs $I_1,...,I_k \subseteq I$ with $I_i \cap A = \emptyset$ for all $i$ and $\sum_{i=1}^k |I_i| \geq \varepsilon |I|$. Finally, if no such $k$ exists, we let $F_{\delta,\varepsilon,A}(n) = \infty$. Note that either there exists an $n \in \bN$ such that $F_{\delta,\varepsilon,A}(m) = 0$ for all $m \geq n$ or else $F_{\delta,\varepsilon,A}$ is an increasing function. 

\begin{thm}
$A$ has the SBM property if and only if: for all $\varepsilon > 0$, there is $\delta > 0$ such that $\lim_{n \to \infty} F_{\delta,\varepsilon,A}(n) = \infty$. 
\end{thm}

We now turn to a discussion of examples and non-examples of sets with the SBM property.  First, an easy example:  

\begin{example}
If $A \subseteq G$ is syndetic, then $A$ has the SBM property. 
\end{example}

\begin{proof}
Assume $A \subseteq G$ is syndetic. Then there is $n \in \bN$ such that $B(g,n) \cap A \neq \emptyset$ for all $g \in G$. By transfer, this will also be true of $^*A$ for all $g \in $ $^*G$. Thus, $g_A(I) \approx 0$ and $\lambda_I(A) = 1$ for all infinite, hyperfinite balls $I \subseteq$ $^*G$, whence $A$ is SBM. 
\end{proof}

In order to proceed further and to see that the SBM property behaves in a way that mimics the SIM property, we will ask that our groups satisfy one further property.

\begin{defn}
We say that the finitely generated group $G$ has the \textbf{small gap property} (or \textbf{SG property} for short) if, for any $r \in (0,1)$, there is $\delta > 0$, such that $\lim_{n \to \infty} \eta_n(G,\delta)>r,$
where $$\eta_n(G,\delta)=\min\{\rho(A,B) \ : \ g_B(A) < \delta, A,B \text{ subballs of } G, |A| > n\}$$ and $$\rho(A,B)=\left\{\frac{|C|}{|A|} : C \subseteq A \cap B \text{ is a subball}\right\}.$$ 
\end{defn}

In other words, if $G$ has the SG property and we know that the gap size between two balls $A$ and $B$ is small, then as long as $A$ is large enough, this also tells us that there is a ball in the intersection that is close in size to $A$. 

 It is unclear at this time what groups satisfy this property and whether or not this property is dependent on the generating set. Certainly, $\bZ^d$ will have this property with its standard generators, but we have yet to verify any non-abelian examples. That being said, it seems possible that the following holds.

\begin{conj}
\label{NiceIntersect}
Let $G$ be a finitely generated, virtually nilpotent group with the SS property. Then $G$ has the SG property. 
\end{conj}

From this point forward, besides assuming that our groups are virtually nilpotent and have the SS property, \textbf{we also assume that our groups have the SG property}, and will make clear when this is used. One consequence of the SG property is that an infinite ball will always have the BM property on any other infinite ball, formalized in the following lemma.

\begin{lem}
If $A$ and $B$ are infinite, hyperfinite balls in $^*G$, then $B$ has the BM property on $A$. 
\end{lem}

\begin{proof}
Fix $\varepsilon \in (0,1)$ and let $\delta>0$ be as in the definition of the SG property for $r = 1-\varepsilon$. Assume that the infinite hyperfinite subball $A' \subseteq A$ is such that $g_B(A') < \delta$. Then the SG property and transfer tell us that there is a ball $C \subseteq A' \cap B$ with $\operatorname{st}\big(\frac{|C|}{|A'|}\big) \geq 1-\epsilon$. By Corollary \ref{NiceMeasures}, 
$$
\lambda_{A'}(B) \geq \lambda_{A'}(C) = \operatorname{st}\left(\frac{|C|}{|A'|}\right) > 1-\varepsilon.
$$ 
Thus, $B$ has the SBM property on $A$. 
\end{proof}

The following lemma is another consequence of the SG property. First, we need to recall the generalization of piecewise syndeticity to arbitrary groups:  $B \subseteq G$ is \textbf{piecewise syndetic} if there is a finite set $F \subseteq G$ such that, for all finite $T \subseteq G$, $gT \subseteq BF$ for some $g \in G$.

\begin{thm}
\label{PWS}
If $B$ is piecewise syndetic, then there is $A \subseteq B$ with the SBM property. 
\end{thm}

We defer the proof of Theorem \ref{PWS} for later, as it will follow immediately from Theorem \ref{BigProof} below. Theorem \ref{PWS} implies that many sets contain SBM sets. However, containing an SBM set is not enough to say that a set is SBM. In fact, the following holds.

\begin{lem}
Let $A \subseteq G$ be a non-syndetic set. Then, there is $B \supseteq A$ such that $B$ is not SBM. 
\end{lem}

\begin{proof}
For every $n \in \bN$, let $x_n \in G$ be such that $B(x_n,n^2) \cap A = \emptyset$. Then set 
$$
B = A \cup \bigcup_{n \in \bN} \bigcup_{0 \leq k \leq n} S(x_n,kn).
$$

Assume that $c$ is a constant witnessing the polynomial growth of $G$ of degree $d$. Fix an infinite $N \in$ $^*\bN$. Fix $\delta > 0$ and take $m \in \bN$ such that $m^d > \frac{c^2}{\delta}$. Let $I = B(x_N,N^2)$, and consider $J = B(x_N,mN)$. By the definition of $B$, we have
$$
g_{^*B}(J) \leq \frac{|B(x_N,N)|}{|B(x_N,mN)|} \leq \frac{cN^d}{\frac{1}{c}m^dN^d} = \frac{c^2}{m^d} < \delta.
$$ 
On the other hand, $^*B \cap J = \bigcup_{0 \leq k \leq m} S(x_N,kN)$ is a finite union of $\lambda_{x_N,mN}$-measure 0 sets by the SS property. Thus, for every $\delta>0$, we can find a $J \subseteq I$ with $g_{^*B}(J) < \delta$ and $\lambda_J($ $^*B) = 0$. Therefore $^*B$ is not BM on $I$ and thus $B$ is not SBM. 
\end{proof}

Note that this lemma only made use of the SS property, and did not require the SG property. However, assuming we do have the SG property, we can combine the previous lemma with Theorem \ref{PWS} to conclude that every set that is piecewise syndetic but not syndetic contains an SBM set, but is also contained in a set that is not SBM. 

We end this subsection with one other fact about SBM sets.

\begin{cor}
If $A$ has the SBM property, then there is an $r \in \bN$ such that for all $x \in Z(G)$, $B(x,r) \cap D(A) \neq \emptyset$.
\end{cor}

\begin{proof}
This follows immediately from Corollary \ref{BMD(A)} and the fact that $D(A) = D($ $^*A)$. 
\end{proof}

\subsection{Supra-SBM Sets}

We stated in the previous section that all piecewise syndetic sets contain an SBM set, and that SBM sets are not closed under supersets. The following definition is thus natural:

\begin{defn}
A set $A \subseteq G$ is \textbf{supra-SBM} if there is a $B \subseteq A$ such that $B$ has the SBM property.
\end{defn}

 With this definition (and making use of the SG property), we prove the following, generalizing a result of Goldbring and Leth \cite[Theorem 3.4]{LGsupraSIM}.

\begin{thm}
\label{BigProof}
Let $A \subseteq G$ be such that $^*A$ has the enhanced BM property on some ball $I$. Then $A$ is supra-SBM.

\end{thm}

\begin{proof}
Let $d$ be the degree of polynomial growth of $G$ and let $c$ be the constant witnessing this growth. Fix $k \in \bN$. Let $\delta(k)$ be as in the definition of the small gap property for $r = \min\big\{\frac{1}{5^dc^2},(1-\frac{1}{k})\big\}$. Fix $\delta_k < \min\big(\delta($ $^*A,I,\frac{1}{k}),\delta(k)\big)$. For any $k,n \in \bN$, if we take $J \subseteq I$ an infinite subball such that $g_{^*A}(J) < \delta_k$, then $\lambda_J($ $^*A)\geq 1 - \frac{1}{k}$ by the definition of $\delta_k$. Therefore, if we take $n$ disjoint subballs $J_1,...,J_n$ of $J$ with $J_i\text{ } \cap $ $^*A = \emptyset$ for each $i$, we have $\sum_{i=1}^n \frac{|J_i|}{|J|}<\frac{1}{k}$. In other words, the sum of $n$ disjoint gaps of $^*A$ in $J$ have sizes adding to less than $\frac{|J|}{k}$. This statement is internal and is true for all infinite $J$, so by underflow, there must be some $M_{n,k} \in \bN$ such that whenever $|J| > M_{n,k}$, then $n$ disjoint gaps of $^*A$ on $J$ add up to less than $\frac{|J|}{k}$. Also, since $\lambda_I($ $^*A) > 0$, for every $n \in \bN$ there exists an infinite subball $J$ of $I$ with $g_{^*A}(J) < \frac{1}{n}$. 

Since there is an infinite ball $I$ with all of these internal properties on $^*A$, by transfer, we can define a sequence of pairwise disjoint balls $I_n$ in $G$ with the following properties:

\begin{itemize}

\item Letting $B(e,r_n)$ be the smallest ball around the identiy which contains $I_1,...,I_{n-1}$, we have that $I_n \cap B(e,nr_n) = \emptyset$.

\item $I_n$ has a subball $J$ with radius at least $n$ such that $g_A(J) < \frac{1}{n}$.

\item For all $k \leq n$ and all $J \subseteq I_n$, if $m < n$, $|J| > M_{m,k}$, and $g_A(J) < \delta_k$, then at least $m+1$ disjoint gaps of $A$ on $J$ are necessary for the size of the gaps to add up to $\frac{|J|}{k}$. 

\end{itemize}

Let $B = \bigcup_{n \in \bN} (A \cap I_n)$. We show that $B$ has the SBM property.

Let $H$ be an infinite, hyperfinite subball of $^*G$. We show that $^*B$ has the BM property on $H$ as witnessed by the function $\delta_H(\varepsilon) = \frac{1}{5^dc^2}\delta_k$ for any $k>\frac{2}{\epsilon}$.  To see this, take $J \subseteq H$ such that $g_{^*B}(J) < \frac{1}{5^dc^2}\delta_k$, and assume that $I_M$ is the ball with largest index which intersects $J$. Note that every $I_K$ with $K < M$ satisfies $I_K \subseteq B(e,r_K)$, whereas every element of $I_M$ is outside of $B(e,Kr_K)$. Therefore, the radius of $J$ is at least $\frac{(K-1)r_K}{2}$, and so $\frac{|B(e,r_K)|}{|J|}$ is infinitesimal as $\frac{2r_K}{(K-1)r_K}$ is infinitesimal. Thus, all elements of $I_K$ with $K < M$ are in the equivalence class of the identity modulo $\sim_M$. 
	
	 As constructed, $^*B \cap H \subseteq B(e,r_M) \cup I_M$. We know that $\frac{|B(e,r_M)|}{|J|}$ is infinitesimal and that $g_{B(e,r_M) \cup I_M}(J) \leq g_{^*B}(J) \leq \frac{1}{5^dc^2}\delta_k$. Assume, towards a contradiction, that $g_{I_M}(J) > \delta_k$. Then there is a closed ball $C \subseteq J$ such that $C \cap I_M = \emptyset$ and $|C| > \delta_k|J|$. Let $r_C$ be the radius of $C$ and $x$ be the center of $C$. Take $y \in C$ such that $d(x,y) = r_C$. Let $x'$ and $y'$ be such that $d(x,x') = d(y,y') =  \lfloor \frac{r_C}{4} \rfloor$ and $d(x,y') = d(y,x') = r_C - \lfloor \frac{r_C}{4} \rfloor$. In other words, let $x'$ and $y'$ be points along some path of length $r_C$ from $x$ to $y$ where $x'$ is approximately $\frac{1}{4}$ of the way from $x$ to $y$ and $y'$ is approximately $\frac{3}{4}$ of the way from $x$ to $y$. Since $d(x,y) = r_C$, we must have that $d(x',y') \geq \frac{r_C}{2}$. This tells us that the balls $B(x',\frac{r_C}{5})$ and $B(y',\frac{r_C}{5})$ are disjoint and contained in $C$. Indeed, if $a$ is in the ball around $x'$ and $b$ is in the ball around $y'$, then $d(a,b) \geq \frac{r_C}{10}$. Finally, note that if $B(e,r_M) \cap B(x',\frac{r_C}{5}) \neq \emptyset$, then $B(e,r_M) \cap B(y',\frac{r_C}{5}) = \emptyset$. Thus 
$$
g_{B(e,r_M) \cup I_M}(J) \geq \frac{1}{c|J|}\left(\frac{r_C}{5}\right)^d \geq \frac{|C|}{|J|c^25^d} > \frac{\delta_k }{c^25^d},
$$
contradicting our original inequality. Thus, we must have that $g_{I_M}(J) \leq \delta_k$. 

Using the SG property and our definition of $\delta_k$, we know that there is a ball $J' \subseteq I_M \cap J$ such that $\frac{|J'|}{|J|} > \frac{1}{5^dc^2}$. This gives us that 
$$
g_{^*B}(J') \leq g_{^*B}(J)\left( \frac{|J|}{|J'|}\right) \leq \left(\frac{1}{5^dc^2}\delta_k\right)(5^dc^2) = \delta_k. 
$$

So $g_{^*A}(J') = g_{^*B}(J') \leq \delta_k$. Since $|J'|$ is infinite, it is larger than $M_{m,k}$ for all $m \in \bN$. Therefore, the cardinalities of $m$ disjoint gaps of $^*A$ (and therefore $^*B$) on $J'$ can not add to $\frac{|J'|}{k}$ for any $m \in \bN$. 

Assume now, towards a contradiction, that $\lambda_{J'}(^*B) < 1 - \frac{1}{k}$. Let $J' = B(g',N')$, so that $\lambda_{J'}(^*B) = \lambda_{g',N'}\big(\operatorname{eq}_{N'}(^*B \cap J')\big)$. We know that $\operatorname{eq}_{N'}(^*B \cap J')$ is a closed set, so that $\overline{J_c} = \operatorname{eq}_{N'}(J') \backslash \operatorname{eq}_{N'}(^*B \cap J')$ is an open set with $\lambda_{g',N'}(\overline{J_c}) > \frac{1}{k}$. Since $\overline{J_c}$ is open, for every $\overline{x} \in \overline{J_c}$, we can find an $n_x \in \bN$ such that $B_{N'}(\overline{x},\frac{1}{n_x}) \subseteq \overline{J_c}$. Let $\cB = \bigcup_{\overline{x} \in \overline{J_c}} \{B_{N'}(\overline{x},\frac{1}{n}) : n \geq n_x\}$. By the symmetric Vitali property, there is a countable pairwise disjoint subcollection $\cB' = \big\{\overline{B_n} : n \in \bN\big\}$ of $\cB$ such that $\lambda_{g',N'}\big(\overline{J_c} \backslash \bigcup_{n \in \bN} \overline{B_n}\big) = 0$, whence $\lambda_{g',N'} \big(\bigcup_{n \in \bN} \overline{B_n}\big) = \lambda_{g',N'}(\overline{J_c}) > \frac{1}{k}$. By countable additivity, we have that $\sum_{n=1}^{\infty} \lambda_{g',N'} (\overline{B_n}) > \frac{1}{k}$. Thus, there must be some $m \in \bN$ such that $\sum_{n=1}^m \lambda_{g',N'} (\overline{B_n}) > \frac{1}{k}$. Because the $\overline{B_n}$ are closed subballs of the open set $\overline{J_c}$ and $\overline{J_c} \cap \operatorname{eq}_{N'}($ $^*B) = \emptyset$, we know that $\operatorname{eq}_{N'}^{-1}(\overline{B_n}) \cap $ $^*B = \emptyset$. Consequently, by overflow, for each $n\leq m$ we can find a closed ball $C_n \subseteq J'$ such that $\operatorname{eq}_{N'}^{-1}(\overline{B_n}) \subseteq C_n$ and $C_n \cap \text{}^*B = \emptyset$. This gives us $m$ gaps of $^*B$ in $J'$ such that

$$
\frac{|\bigcup_{n \leq m} C_n|}{|J'|} \geq \sum_{n \leq m} \lambda_{g',N'}(\overline{B_n}) > \frac{1}{k}.
$$

This contradicts the fact that no finite number of gaps of $^*B$ can add up to $\frac{|J'|}{k}$. Thus, $\lambda_{J'}($ $^*B) \geq 1 - \frac{1}{k}$. 

Finally, $\lambda_J($ $^*B) \geq \lambda_J($ $^*B \cap J') = \operatorname{st}\left(\frac{|J'|}{|J|}\right)\lambda_{J'}($ $^*B) \geq (1-\frac{1}{k})^2$. We conclude that $\lambda_{J}(^*B) > 1-\frac{2}{k}$, as desired. 

It remains to show that $^*B$ has the enhanced BM property on some infinite ball. To see this, note that, for any $N \in $ $^*\bN \backslash \bN$, our construction yields that $I_N$ has a subball $J$ with radius at least $N$  and $g_{^*A}(J) \leq \delta_N < \frac{1}{N} \approx 0$. Since $g_{^*B}(J) = g_{^*A}(J)$, $\lambda_J($ $^*B) = 1$, and so $^*B$ has the enhanced BM property on $J$. 
\end{proof}

Using this theorem and the definition of SBM set, we are able to reformulate the definition of supra-SBM sets in the following way.

\begin{cor}
$A \subseteq G$ is supra-SBM if and only if there is a $B \subseteq A$ and an infinite, hyperfinite ball $I \subseteq $ $^*G$ such that $^*B$ has the enhanced BM property on $I$.
\end{cor}

\begin{proof}
First, assume $A$ is supra-SBM. By definition, there is $B \subseteq A$ with $B$ SBM. This $B$ is as desired by the definition of the SBM property. For the other direction, if $B \subseteq A$ and $^*B$ has the enhanced BM property on $I$, then by Theorem \ref{BigProof} $B$ is supra-SBM, and thus so is $A$. 
\end{proof}

With Theorem \ref{BigProof}, the proof that all piecewise syndetic sets are supra-SBM becomes apparent: 

\begin{proof}[Proof of Theorem \ref{PWS}]
Assume $A \subseteq G$ is piecewise syndetic. Then there is a $d \in \bN$ such that, for every $n \in \bN$, there exists a $g \in G$ with $A$ having gaps of size at most $d$ on $B(g,n)$. Thus, by transfer, there is an infinite, hyperfinite ball $I \subseteq $ $^*G$ such that $^*A$ has gaps of size at most $d$ on $I$. Therefore, for every $J \subseteq I$, $g_{^*A}(J) \approx 0$ and $\lambda_J(^*A) = 1$. By Theorem \ref{BigProof}, $^*A$ has the enhanced BM property on this $I$ and thus $A$ is supra-SBM. 
\end{proof}

Theorem \ref{BigProof} also has the following nice Ramsey-theoretic consequence:

\begin{cor}
Being supra-SBM is a partition regular property of subsets of $G$.
\end{cor}

\begin{proof}
Let $A \subseteq G$ be a supra-SBM set and suppose $A = B_1 \cup \cdots \cup B_n$. Then there is a subset $B$ of $A$ that is SBM and $B = (B \cap B_1) \cup \cdots \cup (B \cap B_n)$. Thus $^*B = $ $^*(B \cap B_1) \cup \cdots \cup $ $^*(B \cap B_n)$, and $^*B$ has the enhanced BM property on some infinite hyperfinite ball $I$. By Theorem \ref{PRBM}, there is an $i\in \{1,\ldots,n\}$ and an infinite subball $J \subseteq I$ such that $^*(B \cap B_i) \cap J$ has the enhanced BM property on $J$. Theorem \ref{BigProof} implies that $B \cap B_i$ is supra-SBM, whence so is $B_i$.
\end{proof}

The last result about supra-SBM sets which we will generalize from the theory of supra-SIM sets is Nathanson's theorem, as discussed in the introduction. Proving this is achieved through repeated use of the following lemma.

\begin{lem}
\label{Kazhdan}
Suppose that $B \subseteq G$ is supra-SBM. Then there is an $r \in \bN$ such that, for every $t \in Z(G)$, there is an element $t' \in B(t,r)$ for which $B \cap Bt'$ is supra-SBM. 
\end{lem}

\begin{proof}
Assume that $A \subseteq B$ is SBM and $^*A$ has the enhanced BM property on the infinite, hyperfinite ball $I$. Let $w \in \bN$ be as in Theorem \ref{DeltaConfigs} for $A_1 = A_2 = $ $^*A$ and $I_1 = I_2 = I$ (for some $\epsilon > 0$ and the corresponding $\delta$). As a reminder, this implies that whenever $J_i = B(a_i,N) \subseteq I$, $i = 1,2$, are such that $g_{^*A}(J_i) \leq \delta$, there is some $c \in G$ such that $^*A \cap B(a_ic,w)\not=\emptyset$. 

Now, fix $t \in Z(G)$ and set 
$$
B_t = \bigcup_{k \in B(e,2w)}t\text{}^*Ak.
$$

\noindent \textbf{Claim:}
If $J$ is an infinite, hyperfinite subball of $I$ such that $\lambda_J(^*A) > 0$, then
$$
^*A \cap B_t \cap J \neq \emptyset.
$$

\textbf{Proof of the Claim:} By the Lebesgue density theorem, there is $a_1 \in J$ and $R \in $ $^*\bN$ such that $\lambda_{B(a_1,R)}(^*A) > 1 - \delta$. By Lemma \ref{ReverseBM}, $g_{^*A}(B(a_1,R)) < \delta$. Note that because $d(a_1,a_1t^{-1}) = |t^{-1}|$, the symmetric difference $B(a_1,R) \triangle B(a_1t^{-1},R)$ is contained in a union of spheres:
$$
B(a_1,R) \triangle B(a_1t^{-1},R) \subseteq \bigcup_{i=0}^{|t^{-1}|} S(a_1,R-i) \cup S(a_1t^{-1},R-i).
$$
Since $t^{-1} \in G$, by the SS property, this is a finite union of measure 0 sets. Thus, we can restrict $^*A$ to $B(a_1,R) \cap B(a_1t^{-1},R)$ without changing our measure. In other words, 
$$
\lambda_{B(a_1,R)}(^*A) = \lambda_{B(a_1,R)}(^*A \cap B(a_1,R) \cap B(a_1t^{-1}R)) =
$$
$$
= \lambda_{B(a_1t^{-1},R)}(^*A \cap B(a_1,R) \cap B(a_1t^{-1}R)) = \lambda_{B(a_1t^{-1},R)}(^*A).
$$

Furthermore, since $t^{-1} \in Z(G)$, we have that 
$$
\lambda_{B(t^{-1}a_1,R)}(^*A) = \lambda_{B(a_1t^{-1},R)}(^*A) = \lambda_{B(a_1,R)}(^*A) > 1-\delta.
$$
 Thus $g_{^*A}(B(t^{-1}a_1,R)) < \delta$. Therefore, we can apply Theorem \ref{DeltaConfigs} to $J_1 = B(a_1,R)$ and $J_2 = B(t^{-1}a_1,R)$, yielding $c \in G$ such that $^*A \cap B(a_1c,w)\not=\emptyset$ and $^*A \cap B(t^{-1}a_1c,w)\not=\emptyset$. Note that 
$$
^*A \cap B(t^{-1}a_1c,w) \neq \emptyset \iff (t\text{ }^*A) \cap B(a_1c,w) \neq \emptyset.
$$
Let $d \in $ $^*A \cap B(a_1c,w) \subseteq $ $^*A \cap J$. Then note that $B(a_1c,w) \subseteq B(d,2w)$ and therefore $B(d,2w) = dB(e,2w)$ contains an element $x$ of $t$ $^*A$. Thus, there is $k^{-1} \in B(e,2w)$ such that $dk^{-1} = x$. This $k$ is also contained in $B(e,2w)$, and so $d \in t$ $^*Ak \subseteq B_t$, finishing the proof of the claim.

The claim implies that, for any infinite hyperfinite subball $J \subseteq I$, we have
$$
\lambda_J(^*A \cap B_t) = \lambda_J(^*A),
$$

as every infinite subball of $J$ either satisfies $\lambda_J(^*A) = 0$ or contains an element of $^*A \cap B_t$. Thus, $^*A \cap B_t$ has the enhanced BM property on $I$. Since $B_t$ is a finite union of internal sets, Theorem \ref{PRBM} implies that there is $k \in B(e,2w)$ such that $^*A \cap t\text{ }^*A k$ has the enhanced BM property on some infinite subball of $I$. It follows that $A \cap tAk$ is supra-SBM, whence, so is $B \cap tBk$. Since $t \in Z(G)$, $B \cap tBk = B \cap Btk$. Thus, for any $t \in Z(G)$, there is $t' = tk \in B(t,2w)$ such that $B \cap Bt'$ is supra-SBM. 
\end{proof}

With this result in place, we prove Nathanson's theorem for supra-SBM sets under the additional assumption that the group has infinite center.

\begin{thm}
\label{NathanSBM}
If $Z(G)$ is infinite and $A \subseteq G$ is supra-SBM, then for every $n \in \bN$, there are $B,C \subseteq G$ such that $B$ is infinite, $|C| = n$, and $BC \subseteq A$. 
\end{thm}

\begin{proof}
First, note that since $Z(G)$ is infinite, for any supra-SBM set $X$, Lemma \ref{Kazhdan} implies there must be infinitely many $t$ such that $X \cap Xt$ is supra-SBM. For $n = 1$, apply Lemma \ref{Kazhdan} to $A$, yielding $t_1 \in G$ such that $B_1=A \cap At_1$ is supra-SBM, and, in particular, is infinite. Letting $C_1 = \{t_1^{-1}\}$, we have our result.
Now, arguing inductively, fix $m \in \bN$ and assume we have constructed $B_m, C_m \subseteq G$ such that $B_m \subseteq A$ is supra-SBM, $C_m = \{t_1^{-1},...,t_m^{-1}\}$, and $B_mC_m \subseteq A$. Then, applying Lemma \ref{Kazhdan} to $B_m$, we find $t_{m+1}\in G$ such that $B_{m+1} = B_m \cap B_mt_{m+1}$ is supra-SBM and $t_{m+1}^{-1} \not\in C_m$. Letting $C_{m+1} = C_m \cup \{t_{m+1}^{-1}\}$, we see that this $B_{m+1}$ and $C_{m+1}$ are as desired.   
\end{proof}

It is shown in \cite{nilpbook} that all finitely generated nilpotent groups have infinite center. In particular, since Theorem \ref{NathanSBM} uses the SS property \textbf{but not} the SG property, it applies to finitely generated groups of nilpotency class 2, such as the discrete Heisenberg group. 

\subsection{Musings on dependence on generators}
As a closing remark, it behooves us to make at least some mention of the role of generating set for the concepts defined in this paper. At the present time, we have been unable to verify whether or not these concepts are independent of generating set. One na\"ive way around this would be to redefine our properties in the following way.

\begin{defn}
We say that a set $A \subseteq G$ has the \textbf{weak SBM (WSBM) property} if there exists a finite, symmetric generating set $S$ of $G$ such that $G$ has the SS and SG properties with respect to this generating set and $A$ is SBM with respect to $S$. Similary, we say that a set $A \subseteq G$ has the \textbf{supra-WSBM property} if there is $B \subseteq A$ such that $B$ has the WSBM property.
\end{defn}

While these definitions are not very illuminating, they do provide us with a generator independent class of sets that satisfy all of the important properties of SBM and supra-SBM sets respectively in virtually nilpotent groups satisfying the SS and SG properties with respect to some generating set. 

It is natural to ask: what sets are or are not WSBM or supra-WSBM? In order to partially answer this question, we need the following:

\begin{lem}
\label{LastLemma}
Fix an internal set $A \subseteq$ $^*G$. Let $S$ and $S'$ be finite generating sets for the virtually nilpotent group $G$, with corresponding word metrics $d_S$ and $d_{S'}$, which in turn induce measures $\lambda_{g,N}$ and $\lambda_{g,N}'$ on the corresponding $d_S$ and $d_{S'}$ balls respectively. Assume $G$ has the SS property with respect to $S'$. Then the existence of an infinite hyperfinite $d_S$ ball $I$ such that $\lambda_I(A) > 0$ implies the existence of an infinite hyperfinite $d_{S'}$ ball $I'$ such that $\lambda_{I'}'(A) > 0$. 
\end{lem}

\begin{proof}
Let $K \in \bN$ be such that the identity map on $G$ is a $K$-bi-Lipschitz map between $(G,d_S)$ and $(G,d_{S'})$.  Let $d$ be the degree of polynomial growth of $G$. Assume that $I$ is an infinite, hyperfinite $d_S$ subball of $G$ such that $\lambda_I(A) > 0$. By the Lebesgue density theorem, for any $\epsilon > 0$, we can find an infinite $d_S$ subball $J = B(g,N) \subseteq I$ with $\lambda_J(A) > 1 - \epsilon$. We know that $J$ contains the $d_{S'}$ ball $I' = B'(g, \lfloor \frac{N}{K} \rfloor)$. Furthermore, there exist constants $c$ and $c'$ such that $|J| \leq cN^d$ and $|I'| \geq \frac{c'}{K^d}N^d$. Thus, we see that $\lambda_J(I') \geq \frac{c'}{cK}$. 

Let $\epsilon = \frac{c'}{2cK}$. Then, since $\lambda_J(A) > 1 - \epsilon$ and $\lambda_J(I') > \epsilon$, we must have that $\lambda_J(A \cap I') > 0$. To conclude the proof, we show that $\lambda_{I'}'(A) > 0$. Assume, towards a contradiction, that $\lambda_{I'}'(A) = 0$. In other words, 
$$
\lambda_{g,\lfloor\frac{N}{K}\rfloor}'\big(\operatorname{eq}_N(A \cap I')\big) = \mu_{g,\lfloor\frac{N}{K}\rfloor}'\big(\operatorname{eq}_N^{-1}\big(\operatorname{eq}_N(A \cap I')\big) \cap I' \big) = 0.
$$
Using the SS property for $G$ with respect to $S'$, we have that 
$$
\mu_{g,\lfloor \frac{N}{K} \rfloor}'\big(\operatorname{eq}_N^{-1}\big(\operatorname{eq}_N(A \cap I')\big)\big) = \mu_{g,\lfloor \frac{N}{K} \rfloor}'\big(\operatorname{eq}_N^{-1}\big(\operatorname{eq}_N(A \cap I')\big) \cap I'\big) = 0.
$$
 Thus, for any $\delta > 0$, there must be an internal set $B$ such that $\operatorname{eq}_N^{-1}\big(\operatorname{eq}_N(A \cap I')\big) \subseteq B$ and $\operatorname{st}\big(\frac{|B|}{|I'|}\big) < \delta$.  Since $B$ will also contain $\operatorname{eq}_N^{-1}\big(\operatorname{eq}_N(A \cap I')\big) \cap J$, $B$ serves to witness
$$
 \lambda_J(A \cap I') = \mu_{g,N}\big(\operatorname{eq}_N^{-1}\big(\operatorname{eq}_N(A \cap I')\big) \cap J\big)\big) \leq \operatorname{st}\Big(\frac{|B|}{|J|}\Big) \leq \operatorname{st}\Big(\frac{|B|}{|I'|}\Big) < \delta. 
$$
Letting $\delta \to 0$, we see that $\lambda_J(A \cap I') = 0$, a contradiction.
\end{proof}

Lemma \ref{LastLemma} shows that many sets in groups with the SS property with respect to some generating set $S$ are not supra-WSBM.

\begin{cor}
\label{LastCor}
Let $G$ be a virtually nilpotent group and let $S$ be a finite generating set for $G$ such that $G$ has the SS property with respect to $S$.  Let $A \subseteq G$ be such that $\lambda_I(^*A) = 0$ for all infinite hyperfinite $d_S$ balls $I \subseteq $ $^*G$. Then $A$ is not supra-WSBM. 
\end{cor} 

\begin{proof}
Fix a finite generating set $S'$ for $G$ with corresponding measure $\lambda_I'$, and let $A$ be as in the statement of the corollary. The contrapositive of Lemma \ref{LastLemma} gives us that $\lambda_I'(^*A) = 0$ for all infinite hyperfinite $d_{S'}$ balls $I$. Therefore, given $B \subseteq A$, we have that $\lambda_I'(^*B) = 0$ for all infinite hyperfinite $d_{S'}$ balls $I$. Thus $^*B$ does not have the enhanced SBM property on any infinite hyperfinite $d_{S'}$ ball.  Since $S'$ is an arbitrary finite generating set for $G$, this yields the desired conclusion.
\end{proof}

\begin{example}
Let $G$ be a finitely generated virtually nilpotent group, and let $S$ be a finite generating set for $G$ such that $G$ has the SS property with respect to $S$. Let $A = \{a_n\}_{n \in \bN} \subseteq G$ be such that $a_1 \neq e$ and $|a_{n+1}| \geq n|a_n|$. Then $A$ is not supra-SSBM.
\end{example}

\begin{proof}
Let $I = B(g,R) \subseteq$ $^*G$ be an infinite hyperfinite ball. By Corollary \ref{LastCor}, it suffices to show that $\lambda_I(^*A)=0$.  By transfer, $|a_{N+1}| \geq N|a_N|$ for all $N\in {}^*\bb N$. Consequently, for all $M \leq N$, we have that
$$
d_S(a_{N+1},a_M) \geq |a_{N+1}| - |a_M| \geq |a_{N+1}| - |a_N| \geq (N-1)|a_N|.
$$ 
This also implies that $\{a_1,...,a_N\} \subseteq B(e,|a_N|)$. 
Let $a_{M_1},...,a_{M_K}$ be all the elements of $^*A$ contained in $I$, with $M_1 < \cdots < M_K$. If there are no such elements, then clearly $\lambda_I(^*A) = 0$. Indeed, we may assume that $K$ is infinite, as otherwise $\lambda_I(^*A) = 0$ as the $\lambda_{g,N}$-measure of finitely many equivalence classes is always 0. 

If $\operatorname{eq}_R(a_{M_K}) = \operatorname{eq}_R(a_{M_{K-1}})$, then $\dfrac{(M_{k}-2)|a_{M_{K}-1}|}{R} \approx 0$, and thus $\operatorname{eq}_R(a_{M_1}) = \cdots = \operatorname{eq}_R(a_{M_K})$. This again gives us that $\lambda_I(^*A) = 0$. 

Finally, assume that $\operatorname{eq}_R(a_{M_K}) \neq \operatorname{eq}_R(a_{M_{K-1}})$. Since $M_{K-1} \leq M_K -1$, we know that 
$$
d_S(a_{M_K},a_{M_{K-1}}) \geq ((M_K-1)-1)|a_{M_K-1}| \geq (M_K-2)|a_{M_{K-1}}|.
$$ 
Thus, in order for $a_{M_K}$ and $a_{M_{K-1}}$ to both be in $B(g,R)$, we must have that 
 $$
 2R \geq d_S(a_{M_K},g) + d_S(a_{M_{K-1}},g) \geq d_S(a_{M_K},a_{M_{K-1}}) \geq (M_K-2)|a_{M_{K-1}}|.
 $$

 Rearranging this inequality, we find 
$$
\frac{|a_{M_{K-1}}|}{R} \leq \frac{2}{M_K-2} \approx 0.
$$
Thus, $\operatorname{eq}_R(B(e,|a_{M_{K-1}}|))$ is a single equivalence class, meaning once again that $\lambda_I(^*A) = 0$ as the measure of two equivalence classes. 
\end{proof}

\end{document}